\documentclass[aos,citesort,dvips]{arximspdf}
\usepackage{mathbh}

\doi{10.1214/10-AOS834}
\volume{39}
\issue{1}
\pubyear{2011}
\firstpage{474}
\lastpage{513}

\makeatletter

\newcommand{\rmd}{{d}}
\newcommand{\nset}{\mathbb{N}}
\newcommand{\rset}{\mathbb{R}}
\newcommand{\thv}{{\theta^\star}}
\newcommand{\calU}{\mathcal U}

\newcommand{\pstat}{\bar{p}}

\newcommand{\eqdef}{\stackrel{\mathrm{def}}{=}}

\newcommand{\rme}{{e}}

\newcommand{\PP}{\mathbb{P}}
\newcommand{\QQ}{\mathbb{Q}}
\newcommand{\PPs}{\bar{\PP}}
\newcommand{\PE}{\mathbb{E}}
\newcommand{\PEs}{\bar{\PE}}
\newcommand{\PPt}{\tilde{\PP}}

\newcommand{\Cset}{\mathsf{C}}
\newcommand{\Xset}{\mathsf{X}}
\newcommand{\Zset}{\mathsf{Z}}
\newcommand{\Xsigma}{\mathcal{X}}
\newcommand{\Hsigma}{\mathcal{H}}
\newcommand{\Fsigma}{\mathcal{F}}
\newcommand{\Zsigma}{\mathcal{Z}}
\newcommand{\JointKernel}{T}

\newcommand{\Yset}{\mathsf{Y}}
\newcommand{\Ysigma}{\mathcal{Y}}
\newcommand{\Q}{Q}
\newcommand{\pt}{\tilde{p}}

\newcommand{\thetapr}{{\theta'}}
\newcommand{\ceta}{\check{\eta}}
\newcommand{\cX}{\check{X}}
\newcommand{\cPP}{\check{\mathbb{P}}}
\newcommand{\cPE}{\check{\mathbb{E}}}
\newcommand{\csigma}{\check{\sigma}}

\newtheorem{theorem}{Theorem}
\newtheorem{prop}[theorem]{Proposition}
\newtheorem{lem}[theorem]{Lemma}

\newproclaim{defi}{Definition}
\newproclaim{rem}{Remark}

\makeatother

\begin{document}
\begin{frontmatter}

\title{Consistency of the maximum likelihood estimator for general
hidden Markov models\thanksref{T1}}
\runtitle{Consistency of the MLE in hidden Markov models}

\thankstext{T1}{Supported in part by the Grant ANR-07-ROBO-0002-04.}

\begin{aug}
\author[A]{\fnms{Randal} \snm{Douc}\ead[label=e1]{randal.douc@it-sudparis.eu}},
\author[B]{\fnms{Eric} \snm{Moulines}\ead[label=e2]{eric.moulines@telecom-paristech.fr}},
\author[C]{\fnms{Jimmy} \snm{Olsson}\ead[label=e3]{jimmy@maths.lth.se}}\\ and
\author[D]{\fnms{Ramon} \snm{van Handel}\corref{}\ead[label=e4]{rvan@princeton.edu}}
\runauthor{Douc, Moulines, Olsson and van Handel}
\affiliation{T{\'e}l{\'e}com SudParis, T{\'e}l{\'e}com ParisTech,
Lund University and Princeton University}
\address[A]{R. Douc\\
CITI/T{\'e}l{\'e}com SudParis\\
9 rue Charles Fourier\\
91000 Evry\\
France\\
\printead{e1}}
\address[B]{E. Moulines\\
CNRS/LTCI/T{\'e}l{\'e}com ParisTech\\
46 rue Barrault\\
75013 Paris\\
France\\
\printead{e2}}
\address[C]{J. Olsson\\
Center of Mathematical Sciences\\
Lund University\\
Box 118\\
SE-22100 Lund\\
Sweden\\
\printead{e3}}
\address[D]{R. van Handel\\
Department of Operations Research \\
\quad and Financial Engineering\\
Princeton University\\
Princeton, New Jersey 08544\\
USA\\
\printead{e4}}
\end{aug}

\received{\smonth{12} \syear{2009}}
\revised{\smonth{4} \syear{2010}}

%
\begin{abstract}
Consider a parametrized family of general hidden Markov models, where
both the observed and unobserved components take values in a
complete separable metric space. We prove that the maximum likelihood
estimator (MLE) of the parameter is strongly consistent under a rather
minimal set of assumptions. As special cases of our main result, we
obtain consistency in a large class of nonlinear state space models, as
well as general results on linear Gaussian state space models and finite
state models.

A novel aspect of our approach is an information-theoretic technique for
proving identifiability, which does not require an explicit
representation for the relative entropy rate. Our method of proof
could therefore form a foundation for the investigation of MLE consistency
in more general dependent and non-Markovian time series. Also of
independent interest is a general concentration inequality for
$V$-uniformly ergodic Markov chains.
\end{abstract}

%

\end{frontmatter}

\section{Introduction}
\label{sec:notations-definitions}

A hidden Markov model (HMM) is a bivariate stochastic process $(X_k,
Y_k)_{k\geq0}$, where $(X_k)_{k\geq0}$ is a Markov chain (often
referred to as the state sequence) in a state space $\Xset$ and,
conditionally on $(X_k)_{k \geq0}$, $(Y_k)_{k\geq0}$ is a sequence of
independent random variables in a state space $\Yset$ such that the
conditional distribution of $Y_k$ given the state sequence depends on
$X_k$ only. The key feature of HMM is that the state sequence
$(X_k)_{k\geq0}$ is not observable, so that statistical inference has
to be carried out by means of the observations $(Y_k)_{k\geq0}$ only.
Such problems are far from straightforward due to the fact that the
observation process $(Y_k)_{k\geq0}$ is generally a dependent,
non-Markovian time series [despite that the bivariate process
$(X_k,Y_k)_{k\geq0}$ is itself Markovian]. HMM appear in a large
variety of scientific disciplines including financial econometrics
\cite{hullwhite1987,mamonelliott2007}, biology
\cite{churchill1992}, speech recognition \cite{juangrabiner1991},
neurophysiology \cite{fredkinrice1987}, etc., and the statistical
inference for such models is therefore of significant practical
importance \cite{cappemoulinesryden2005}.

In this paper, we will consider a parametrized family of HMM with
parameter space $\Theta$. For each parameter $\theta\in\Theta$, the
dynamics of the HMM is specified by the transition kernel $Q_\theta$ of
the Markov process $(X_k)_{k\ge0}$, and by the conditional distribution
$G_\theta$ of the observation $Y_k$ given the signal $X_k$. For
example, the state and observation sequences may be generated according
to a nonlinear dynamical system (which defines implicitly $Q_\theta$ and
$G_\theta$) of the form
\begin{eqnarray*}
X_k &=& a_\theta(X_{k - 1}, W_k) , \\
Y_k &=& b_\theta(X_k, V_k) ,
\end{eqnarray*}
where $a_\theta$ and $b_\theta$ are (nonlinear) functions and
$(W_k)_{k \geq1}$, $(V_k)_{k \geq0}$ are independent sequences of
i.i.d. random variables which are independent of $X_0$.

Throughout the paper, we fix a distinguished element $\thv\in\Theta$.
We will always presume that the kernel $Q_\thv$ possesses a unique
invariant probability measure $\pi_\thv$, and we denote by $\PPs
_\thv$
and $\PEs_\thv$ the law and associated expectation of the stationary HMM
with parameter $\thv$ (we refer to Section \ref{sec:canonical} for
detailed definitions of these quantities). In the
setting of this paper, we have access to a single observation path of
the process $(Y_k)_{k \ge0}$ sampled from the distribution $\PPs_\thv$.
Thus, $\thv$ is interpreted as the \textit{true} parameter value, which is
not known a priori. Our basic problem is to form a consistent estimate
of $\thv$ on the basis of the observations $(Y_k)_{k \ge0}$ only,
that is,
without access to the hidden process $(X_k)_{k \ge0}$. This will be
accomplished by means of the maximum likelihood method.

The maximum likelihood estimator (MLE) is one of the backbones of
statistics, and common wisdom has it that the MLE should be, except in
``atypical'' cases, consistent in the sense that it converges to the
true parameter value as the number of observations tends to infinity.
The purpose of this paper is to show that this is indeed the case for
HMM under a rather minimal set of assumptions. Our main result
substantially generalizes previously known consistency results for HMM,
and can be applied to many models of practical interest.

\subsection{Previous work}

The study of asymptotic properties of the MLE in HMM was initiated in
the seminal work of Baum and Petrie \cite{baumpetrie1966,petrie1969}
in the 1960s. In these papers, the state space $\Xset$ and the
observation space $\Yset$ were both presumed to be finite sets. More
than two decades later, Leroux \cite{leroux1992} proved consistency for
the case that $\Xset$ is a finite set and $\Yset$ is a general state
space. The consistency of the MLE in more general HMM has subsequently
been investigated in a series of contributions
\cite{leglandmevel2000a,leglandmevel2000,doucmatias2002,doucmoulinesryden2004,genoncatalotlaredo2006}
using a variety of methods. However, all
these results require very restrictive assumptions on the underlying
model, such as uniform positivity of the transition densities, which are
rarely satisfied in applications (particularly in the case of a
noncompact state space $\Xset$). A general consistency result for HMM has
hitherto remained lacking.

Though the consistency results above differ in the details of their
proofs, all proofs have a common thread which serves also as the
starting point for this paper. Let us therefore recall the basic
approach for proving consistency of the MLE. Denote by
$p^\nu(y_0^n;\theta)$ the likelihood of the observations $Y_0^n$ for the
HMM with parameter $\theta\in\Theta$ and initial measure $X_0\sim
\nu$.
The first step of the proof aims to establish that for any
$\theta\in\Theta$, there is a constant $H(\thv,\theta)$ such that
\[
\lim_{n \to\infty} n^{-1} \log p^\nu(Y_{0}^{n};\theta)
= \lim_{n \to\infty} n^{-1} \PEs_\thv[
\log p^\nu(Y_{0}^{n};\theta)]
= H(\thv,\theta) ,\qquad \PPs_{\thv}\mbox{-a.s.}
\]
For $\theta=\thv$, this convergence follows from the generalized
Shannon--Breiman--McMillan theorem \cite{barron1985}, but for
$\theta\ne\thv$ the existence of the limit is far from obvious.
Now set $K(\thv,\theta)=H(\thv,\thv)-H(\thv,\theta)$. Then
$K(\thv,\theta)\ge0$ is the relative entropy rate between the observation
laws of the parameters $\thv$ and $\theta$, respectively. The
second step of the proof aims to establish identifiability, that is,
that $K(\thv,\theta)$ is minimized only at those parameters
$\theta$ that are equivalent to $\thv$ (in the sense that they give rise
to the same stationary observation law). Finally, the third step of
the proof aims to prove that the maximizer of the likelihood
$\theta\mapsto p^\nu(Y_0^n;\theta)$ converges
$\PPs_\thv$-a.s. to the maximizer of $H(\thv,\theta)$, that is,
to the minimizer of $K(\thv,\theta)$. Together, these three
steps imply consistency.

Let us note that one could write the likelihood as
\[
n^{-1}\log p^\nu(Y_0^n;\theta) =
\frac{1}{n}\sum_{k=0}^n \log p^\nu(Y_k|Y_0^{k-1};\theta),
\]
where $p^\nu(Y_k|Y_0^{k-1};\theta)$ denotes the conditional density of
$Y_k$ given $Y_0^{k-1}$ under the parameter $\theta$ (i.e., the
one-step predictor). If the limit of $p^\nu(Y_1|Y_{-n}^{0};\theta)$ as
$n\to\infty$ can be shown to exist $\PPs_\thv$-a.s., existence of the
relative entropy rate follows from the ergodic theorem and yields the
explicit representation $H(\thv,\theta) = \PEs_\thv[\log
p^\nu(Y_1|Y_{-\infty}^{0};\theta)]$. Such an approach was used in
\cite{baumpetrie1966,doucmoulinesryden2004}. Alternatively, the
predictive distribution $p^\nu(Y_k|Y_0^{k-1};\theta)$ can be expressed
in terms of a measure-valued Markov chain (the prediction filter), so
that existence of the relative entropy rate, as well as an explicit
representation for $H(\thv,\theta)$, follows from the ergodic theorem
for Markov chains if the prediction filter can be shown to be ergodic.
This approach was used in \cite{leglandmevel2000a,leglandmevel2000,doucmatias2002}.
In \cite{leroux1992}, the existence of the relative
entropy rate is established by means of Kingman's subadditive ergodic
theorem (the same approach is used indirectly in \cite{petrie1969},
which invokes the Furstenberg--Kesten theory of random matrix products).
After some additional work, an explicit representation of
$H(\thv,\theta)$ is again obtained. However, as noted in \cite{leroux1992},
page 136, the latter is surprisingly difficult, as Kingman's
ergodic theorem does not directly yield a representation of the limit as
an expectation.

Though the proofs use different techniques, all the results above rely
heavily on the explicit representation of $H(\thv,\theta)$ in order to
establish identifiability. This has proven to be one of the main
difficulties in developing consistency results for more general HMM. For
example, an attempt in \cite{genoncatalotlaredo2006} to generalize
the approach of \cite{leroux1992} failed to establish such a
representation, and therefore to establish consistency except in a
special example. Once identifiability has been established, standard
techniques (such as Wald's method) can be used to show
convergence of the maximizer of the likelihood, completing the proof.

For completeness, we note that a recent attempt \cite{fuh2006} to prove
consistency of the MLE for general HMM contains very serious problems in
the proof \cite{jensen2009} (not addressed in \cite{fuh2009}), and
therefore fails to establish the claimed results.

\subsection{Approach of this paper}

In this paper, we prove consistency of the MLE for general HMM under
rather mild assumptions. Though our proof follows broadly the general
approach described above, our approach differs from previous work in two
key aspects. First, we note that it is not necessary to establish
existence of the relative entropy rate. Indeed, rather than attempting
to prove the existence of a limiting contrast function
\[
\lim_{n \to\infty} n^{-1} \log p^\nu(Y_0^n;\theta)
= H(\thv,\theta) , \qquad\PPs_{\thv}\mbox{-a.s.},
\]
which must then shown to be identifiable in the sense that
$H(\thv,\theta)<H(\thv,\thv)$ for parameters $\theta$ not equivalent
to $\thv$, it suffices to show directly that
\[
\limsup_{n \to\infty} n^{-1} \log p^\nu(Y_0^n;\theta)
< H(\thv,\thv) ,\qquad \PPs_{\thv}\mbox{-a.s.}
\]
[note that the existence of $H(\thv,\thv)$ is guaranteed by the
Shannon--Breiman--McMillan theorem, and therefore poses little
difficulty in
the proof]. This simple observation implies that it suffices to obtain
a convenient upper bound for $p^\nu(Y_0^n;\theta)$, which we
accomplish by introducing the assumption that some iterate $Q^l_\theta$
of the transition kernel of the state sequence possesses a bounded
density with respect to a $\sigma$-finite reference measure $\lambda$.

Second, and perhaps more importantly, we avoid entirely the need to
obtain an explicit representation for the limiting contrast function
$H(\thv,\theta)$ which played a key role in all previous work. Instead,
we develop in Section \ref{sec:identifiability} a surprisingly powerful
information-theoretic device which may be used to prove identifiability
in a very general setting (see \cite{martonshields1994} for related
ideas), and is not specific to HMM. This technique yields the
following: in order to establish that the normalized relative entropy is
bounded away from zero, that is,
\[
\liminf_{n \to\infty} \PEs_{\thv} \biggl[ n^{-1}
\log\frac{\pstat(Y_0^n;\thv)}{
p^\nu(Y_0^n;\theta)} \biggr] > 0 ,
\]
it suffices to show that there is a sequence of sets $(A_k)_{k\ge0}$
such that
\[
\liminf_{n\to\infty}\PPs_\thv(Y_0^n\in A_n)>0,\qquad
\limsup_{n\to\infty}n^{-1}\log\PP_\theta^\nu(Y_0^n\in A_n)<0
\]
[here $\PP_\theta^\nu$ is the law of the HMM with parameter $\theta
$ and
initial measure $\nu$, while $\pstat(y_0^n;\thv)$ denotes the likelihood
of $Y_0^n$ under $\PPs_\thv$]. It is rather straightforward to find such
a sequence of sets, provided the law of the observations $(Y_k)_{k\ge0}$
is ergodic under $\PPs_\thv$ and satisfies an elementary large deviations
property under $\PP_\theta^\nu$. These properties are readily established
in a very general setting. In particular, we will show (Section
\ref{sec:v-uniform-separation}) that any geometrically ergodic state
sequence gives rise to the requisite large deviations property, so that
our main result can be applied immediately to a large class of models of
practical interest. (Let us note, however, that ergodicity of
$\PP_\theta^\nu$ is not necessary; see Section \ref
{sec:finite-state-leroux}.)

Of course, there are some complications. Rather than investigating the
likelihood function $p^\nu(Y_0^n;\theta)$ directly, the proof
of our main result relies in an essential manner on the asymptotics of
the process $p^\lambda(Y_0^n;\theta)$ where $\lambda$ is the
reference measure defined above. The latter process plays a special role
in our proofs due to the fact that it satisfies a certain
submultiplicativity property; this allows us to upper bound $n^{-1}\log
p^\nu(Y_0^n;\theta)$ by a time average, which possesses an
almost sure limit by Birkhoff's ergodic theorem (see the proof of
Theorem \ref{thm:consistency} below for further details). As $\lambda$
is typically only $\sigma$-finite, however, it is not immediately
obvious that the problem is well-posed. Nonetheless, we will see that
these complications can be resolved, provided that the HMM is
sufficiently ``observable'' so that the improper likelihood
$p^\lambda(Y_0^n;\theta)$ is well defined for sufficiently
large $n$ (under mild integrability conditions). As is demonstrated by
the examples in Section \ref{sec:applications}, this is the case in a
wide variety of applications.

Finally, let us note that the techniques used in the proof of our main
result appear to be quite general. Though we have restricted our
attention in this paper to the case of HMM, these techniques could form
the foundation for consistency proofs in other dependent and
non-Markovian time series models (such as, e.g., the
autoregressive setting of \cite{doucmoulinesryden2004}), which share
many of the difficulties of statistical inference in hidden Markov
models. Other asymptotic properties of the MLE, such as asymptotic
normality, merit further investigation.

\subsection{Organization of the paper}

The remainder of the paper is organized as follows. In Section
\ref{sec:AssumptionsMainResults}, we first introduce the setting and
notations that are used throughout the paper. Then, we state our main
assumptions and results. In Section \ref{sec:applications}, our main
result is used to establish consistency in three general classes of
models: linear-Gaussian state space models, finite state models, and
nonlinear state space models of the vector ARCH type (this includes the
stochastic volatility model and many other models of interest in time
series analysis and financial econometrics).
Section \ref{sec:main-proof} is devoted to the proof of our main
result. Finally, Section~\ref{sec:v-uniform-separation} is devoted to
the proof of the fact that geometrically ergodic models satisfy the
large deviations property needed for identifiability. In particular, we
prove in Section \ref{sec:exponential-inequality} general
Azuma--Hoeffding type concentration inequality for $V$-uniformly
ergodic Markov chains, which is of independent interest.

\section{Assumptions and main results}
\label{sec:AssumptionsMainResults}

\subsection{Canonical setup and notation}
\label{sec:canonical}

We fix the following spaces throughout:
\begin{itemize}
\item$\Xset$ is a Polish space endowed with
its Borel $\sigma$-field $\Xsigma$.
\item$\Yset$ is a Polish space endowed with
its Borel $\sigma$-field $\Ysigma$.
\item$\Theta$ is a compact metric space endowed with
its Borel $\sigma$-field $\Hsigma$.
\end{itemize}
$\Xset$ is the state space of the hidden Markov process, $\Yset$ is the
state space of the observations, and $\Theta$ is the parameter space of
our model. We furthermore assume that $\Theta$ is endowed with a given
equivalence relation\setcounter{footnote}{1}\footnote{This
is meant here in the broad sense, that is, $\sim$ is a
binary relation on $\Theta$ indicating which elements
$\theta\in\Theta$ should be viewed as ``equivalent.'' We
do not require $\sim$ to be transitive.

It should be emphasized that in the setting of this paper,
the equivalence relation $\sim$ is presumed to be given
as part of the model specification, rather than being defined
in terms of the model: the statistician may choose up to
which equivalence she wishes to estimate the true parameter
generating the observations. One assumption of our main result
[assumption (A6) below] then requires that
parameters $\theta,\thetapr$ that are not equivalent,
denoted $\theta\not\sim\thetapr$, give rise to observation
laws that are distinguishable in a suitable sense.
In many cases, there is a natural equivalence relation
which ensures that this is the case; see Section
\ref{sec:sufficient-exponential-separation} below.} $\sim$, and denote the equivalence class of $\theta\in\Theta$ as
$[\theta]\eqdef\{\theta'\in\Theta\dvtx\theta'\sim\theta\}$.

Our model is defined as follows: we are given a transition kernel
$Q\dvtx\Theta\times\Xset\times\Xsigma\to[0,1]$, a positive $\sigma$-finite
measure $\mu$ on $(\Yset,\Ysigma)$, and a measurable function
$g\dvtx\Theta\times\Xset\times\Yset\to\mathbb{R}_+$ such that
$\int g_\theta(x,y) \mu(\rmd y)=1$ for all $\theta,x$. For each
$\theta\in\Theta$, we can define the transition kernel $T_\theta$
on $(\Xset,\Yset)$ as
\[
\JointKernel_\theta[(x,y), C] \eqdef
\int\mathbh{1}_C(x',y') g_\theta(x',y') \mu(\rmd y')
\Q_\theta(x,\rmd x').
\]
We will work on the measurable space $(\Omega,\Fsigma)$ where
$\Omega=(\Xset\times\Yset)^\nset$,
$\Fsigma=(\Xsigma\otimes\Ysigma)^{\otimes\nset}$, and the canonical
coordinate process is denoted as $(X_k,Y_k)_{k\ge0}$. For each
$\theta\in\Theta$ and probability measure $\nu$ on $(\Xset,\Xsigma
)$, we
define $\PP^\nu_\theta$ to be the probability measure on
$(\Omega,\Fsigma)$ such that $(X_k,Y_k)_{k\ge0}$ is a time homogeneous
Markov process with initial measure $\PP^\nu_\theta((X_0,Y_0)\in
C)=\int
\mathbh{1}_C(x,y) g_\theta(x,y)\* \mu(\rmd y) \nu(\rmd x)$ and transition
kernel $T_\theta$. Denote as $\PE^\nu_\theta$ the expectation with
respect to $\PP^\nu_\theta$, and denote as $\PP^{\nu,Y}_\theta$ the
marginal of the probability measure $\PP^\nu_\theta$ on
$(\Yset^\nset,\Ysigma^{\otimes\nset})$.

Throughout the paper, we fix a distinguished element $\thv\in\Theta
$. We
will always presume that the kernel $Q_\thv$ possesses a unique
invariant probability measure $\pi_\thv$ on $(\Xset,\Xsigma)$ [this
follows from assumption (A1) below]. For ease of
notation, we will write $\PPs_\thv,\PEs_\thv,\PPs_\thv^Y$ instead of
$\PP_\thv^{\pi_\thv},\PE_\thv^{\pi_\thv},\PP_\thv^{\pi_\thv,Y}$. Though
the kernel, $Q_\theta$ need not be uniquely ergodic for $\theta\ne
\thv$
in our main result, we will obtain easily verifiable assumptions in a
setting which implies that all $Q_\theta$ possess a unique invariant
probability measure. When this is the case, we will denote as
$\pi_\theta$ this invariant measure and we define
$\PPs_\theta,\PEs_\theta,\PPs_\theta^Y$ as above.

Under the measure $\PP^\nu_\theta$, the process $(X_k,Y_k)_{k\ge0}$ is
a \textit{hidden Markov model}. The hidden process $(X_k)_{k\ge0}$ is a
Markov chain in its own right with initial measure $\nu$ and transition
kernel $Q_\theta$, while the observations $(Y_k)_{k\ge0}$ are
conditionally independent given the hidden process with common
observation kernel $G_\theta(x,\rmd y)=g_\theta(x,y) \mu(\rmd y)$. In
the setting of this paper, we have access to a single observation path
of the process $(Y_k)_{k\ge0}$ sampled from the distribution $\PPs
_\thv$.
Thus, $\thv$ is interpreted as the \textit{true} parameter value, which is
not known a priori. Our basic problem is to obtain a consistent estimate
of $\thv$ (up to equivalence, i.e., we aim to identify the equivalence
class $[\thv]$ of the true parameter) on the basis of the observations
$(Y_k)_{k\ge0}$ only, without access to the hidden process
$(X_k)_{k\ge0}$.
This will be accomplished by the maximum likelihood method.

Define for any positive $\sigma$-finite measure $\rho$ on $(\Xset
,\Xsigma)$
\begin{eqnarray*}
p^\rho(\rmd x_{t+1},y_s^t; \theta)
&\eqdef&\int\rho(\rmd x_s) \prod_{u = s}^t
g_{\theta}(x_u, y_u)
Q_\theta(x_u, \rmd x_{u+1})
, \\
p^\rho(y_s^t; \theta) &\eqdef&
\int p^\rho(\rmd x_{t+1},y_s^t; \theta),
\end{eqnarray*}
with the conventions $\prod_{u=v}^w a_u=1$ if $v>w$ and for any sequence
$(a_s)_{s \in\mathbb{Z}}$ and any integers $s\leq t$, $a_s^t\eqdef
(a_s, \ldots, a_t)$. For ease of notation, we will write
$p^x(y_s^t; \theta)\eqdef p^{\delta_x}(y_s^t;\theta)$
for $x\in\Xset$, and we write
$\pstat(y_s^t; \theta) \eqdef
p^{\pi_\theta}(y_s^t;\theta)$. Note that
$p^\rho(\rmd x_{t+1},y_s^t; \theta)$ is a positive but not
necessarily $\sigma$-finite measure. However, if $\rho$ is a probability
measure, then $p^\rho(\rmd x_{t+1},y_s^t; \theta)$ is a finite
measure and $p^\rho(y_s^t;\theta)<\infty$.

If $\nu$ is a probability measure, then $p^\nu(y_0^n;\theta)$ is the
likelihood of the observation sequence $y_0^n$ under the law
$\PP_\theta^\nu$. The maximum likelihood method forms an estimate of
$\thv$ by maximizing $\theta\mapsto p^\nu(y_0^n;\theta)$, and we
aim to
establish consistency of this estimator. However, as the state space
$\Xset$ is not compact, it will turn out to be essential to consider also
$p^\lambda(y_0^n;\theta)$ for a positive $\sigma$-finite measure
$\lambda$.

We conclude this section with some miscellaneous notation. For any
function~$f$, we denote as $|f|_\infty$ its supremum norm [e.g.,
$|g_{\theta}|_\infty\eqdef\sup_{(x,y) \in\Xset\times\Yset}
g_\theta(x$, $y)$]. As we will frequently integrate with respect to the
measure $\mu$, we will use the abridged notation $\rmd y$ instead of
$\mu(\rmd y)$, and we write $\rmd y_s^t\eqdef\prod_{i=s}^t
\rmd y_i$. For any integer $m$ and $\theta\in\Theta$, we denote by
$Q_\theta^m$ the $m$th iterate of the kernel $Q_\theta$. For any pair
of probability measures $\PP,\QQ$ and function $V\ge1$, we define the
norm
\[
\|\PP-\QQ\|_V \eqdef\sup_{f\dvtx|f|\le V}\biggl|
\int f \,d\PP- \int f \,d\QQ
\biggr|.
\]
Finally, the relative entropy (or Kullback--Leibler divergence) is
defined as
\[
\mathrm{KL}(\PP||\QQ)\eqdef
\cases{\displaystyle \int\log(\rmd\PP/\rmd\QQ) \,\rmd\PP, &\quad if $\PP\ll\QQ$,\vspace*{2pt}\cr
\infty, &\quad otherwise,}
\]
for any pair of probability measures $\PP$ and $\QQ$.
\begin{rem}
Throughout the paper, we will encounter partial suprema of measurable
functions [e.g., $y_0^n\mapsto\sup_{\theta\in\mathcal{U}}
\pstat(y_0^n;\theta)$ for some measurable set $\mathcal{U}\in
\Hsigma$].
As the supremum is taken over an uncountable set, such functions are not
necessarily Borel-measurable. However, as all our state spaces are
Polish, such functions are always guaranteed to be universally
measurable (\cite{bertsekasshreve1978}, Proposition 7.47). Similarly,
a~Borel-measurable (approximate) maximum likelihood estimator need not
exist, but the Polish assumption ensures the existence of universally
measurable maximum likelihood estimators
(\cite{bertsekasshreve1978}, Proposition 7.50). All probabilities and
expectations can therefore be unambiguously extended to such quantities,
which we will implicitly assume to be the case in the sequel.
\end{rem}

\subsection{The consistency theorem}

Our main result establishes consistency of the MLE under assumptions
(A1)--(A6) below, which hold in a
large class of models. Various examples will be treated in Section
\ref{sec:applications} below.

\begin{enumerate}[(A6)]
\item[(A1)] 
The Markov kernel $Q_{\thv}$ is positive Harris recurrent.
\end{enumerate}

\begin{enumerate}[(A6)]
\item[(A2)] 
$\PEs_\thv[\sup_{x\in\Xset} (\log g_\thv(x,Y_0))^+]<\infty$,
$\PEs_\thv[|{\log\int g_\thv(x,Y_0) \pi_\thv(\rmd
x)}|]
<\infty$.
\end{enumerate}

Assumptions (A1), (A2) ensure the
existence of the entropy rate for $\thv$.

\begin{enumerate}[(A6)]
\item[(A3)] 
There is an integer $l\ge1$, a measurable function
$q\dvtx\Theta\times\Xset\times\Xset\to\mathbb{R}_+$, and a $\sigma$-finite
measure $\lambda$ on $(\Xset,\Xsigma)$
such that $|q_{\theta}|_\infty< \infty$ and
\[
Q_\theta^l(x,A) = \int\mathbh{1}_A(x') q_\theta(x,x') \lambda
(\rmd x')
\]
for all $\theta\not\sim\thv$, $x\in\Xset$, $A\in\Xsigma$.
\end{enumerate}

Assumption (A3) states that an iterate of the transition
kernel $Q_\theta$ possesses a density with respect to a $\sigma$-finite
measure $\lambda$. This property will allow us to establish the
asymptotics of the likelihood of $\PP_\theta^\nu$ in terms of the
improper likelihood $p^\lambda(\cdot;\theta)$. The measure $\lambda$
plays a central role throughout the paper.

\begin{enumerate}[(A6)]
\item[(A4)] 
For every $\theta\not\sim\thv$, there is a neighborhood $\calU
_\theta$ of
$\theta$ such that
\[
\sup_{\thetapr\in\calU_\theta}|q_{\thetapr}|_\infty<\infty,\qquad
\PEs_\thv\Bigl[\sup_{\thetapr\in\calU_\theta}
\sup_{x\in\Xset} (\log g_\thetapr(x,Y_0))^+\Bigr]
<\infty,
\]
and there is an integer $r_\theta$ such that
\[
\PEs_\thv\Bigl[\sup_{\thetapr\in\calU_\theta}
(\log p^\lambda({Y}_{0}^{r_\theta}; \thetapr))^+\Bigr]
<\infty.
\]
\end{enumerate}

\begin{enumerate}[(A6)]
\item[(A5)] 
For any $\theta\not\sim\thv$ and $n \geq r_\theta$, the function
$\thetapr\mapsto p^\lambda(Y_0^n;\thetapr)$ is
upper-semicontinuous at $\theta$, $\PPs_\thv\mbox{-a.s.}$
\end{enumerate}

Assumptions (A4) and (A5) are similar in spirit
to the classical Wald conditions in the case of i.i.d. observations.
However, an important difference with the classical case is that
(A4) applies to
$p^\lambda({y}_{0}^{r_\theta};\theta)$, which is not a probability
density (as $\lambda$ is typically only $\sigma$-finite). Assumption
(A4) implies in particular that
$p^\lambda({y}_{0}^{r_\theta};\theta)$ is $\PPs_{\thv}\mbox
{-a.s.}$ finite.
When $\lambda$ is $\sigma$-finite, this requires, in essence, that the
observations contain some information on the range of values taken by
the hidden process.

Finally, the key assumption (A6) below gives identifiability of
the model. In principle, what is needed is that $\PP_\theta^{\lambda
,Y}$ is
distinguishable from $\PPs_\thv^Y$ in a suitable sense. However, as
$\lambda$ may be $\sigma$-finite, $\PP_\theta^{\lambda,Y}$ is not well
defined. As a replacement, we will consider the probability measure
$\PPt^\lambda_\theta$ defined by
%
%
\begin{equation}
\PPt^\lambda_\theta(Y_0^n\in A) = \int\mathbh{1}_A(y_0^n)
\frac{
p^\lambda({y}_{0}^{n};\theta)}{
p^\lambda({y}_{0}^{r_\theta};\theta)}
\pstat({y}_{0}^{r_\theta};\thv)\,
\rmd y_0^n
\end{equation}
for all $n\ge r_\theta$ and $A\in\Ysigma^{\otimes(n+1)}$ (note that the
definition of $\PPt^\lambda_\theta$ depends implicitly on $\thv$ as well
as on $\theta$; the former dependence is suppressed for notational
simplicity). Lemma \ref{lem:ptilde} shows that $\PPt^\lambda_\theta
$ is
well defined, provided that (A4) holds and
$p^\lambda({Y}_{0}^{r_\theta};\theta)>0$ $\PPs_\thv\mbox{-a.s.}$
The law $\PPt^\lambda_\theta$ is in essence a normalized version of
$\PP_\theta^{\lambda,Y}$, and (A6) should be interpreted
in this spirit.

\begin{enumerate}[(A6)]
\item[(A6)] 
For every $\theta\not\sim\thv$ such that
$p^\lambda({Y}_{0}^{r_\theta};\theta)>0$ $\PPs_\thv\mbox{-a.s.}$,
we have
\[
\liminf_{n\to\infty}\PPs_\thv(Y_0^n\in A_n)>0,\qquad
\limsup_{n\to\infty}n^{-1}\log\PPt^\lambda_\theta(Y_0^n\in A_n)<0
\]
for some sequence of sets $A_n\in\Ysigma^{\otimes(n+1)}$.
\end{enumerate}

Despite that this assumption looks nontrivial, we will obtain sufficient
conditions in Section \ref{sec:sufficient-exponential-separation}
which are satisfied in a large class of models.

Having introduced the necessary assumptions, we now turn to the
statement of our main result. Let $\ell_{\nu,n}\dvtx \theta\mapsto\log
p^\nu(Y_0^n; \theta)$ be the log-likelihood function
associated with the initial probability measure $\nu$ and the
observations $Y_0^n$. An approximate \textit{maximum likelihood
estimator} $({\hat{\theta}}_{\nu,n})_{n\ge0}$ is defined as a
sequence of
(universally) measurable functions ${\hat{\theta}}_{\nu,n}$ of
$Y_0^n$ such that
\[
n^{-1}\ell_{\nu,n}({\hat{\theta}}_{\nu,n}) \ge
\sup_{\theta\in\Theta}n^{-1}\ell_{\nu,n}(\theta)
- o_{\mathrm{a.s.}}(1),
\]
where\vspace*{2pt} $o_{\mathrm{a.s.}}(1)$ denotes a stochastic process that converges to
zero $\PPs_\thv$-a.s.\break as $n\to\infty$ [if the supremum of
$\ell_{\nu,n}$ is attained, we may choose
${\hat{\theta}}_{\nu,n}=\break \arg\max_{\theta\in
\Theta}
\ell_{\nu,n}(\theta)$]. The main result of the paper consists
in obtaining the consistency of ${\hat{\theta}}_{\nu,n}$.
\begin{theorem}\label{thm:consistency}
Assume \textup{(A1)--(A6)}, and let $\nu$ be a fixed
initial probability measure. Suppose that one of the following assumptions
hold:
\begin{enumerate}
\item$\nu\sim\pi_\thv$; or
\item$g_\thv(x,y)>0$ for all $x,y$, and
$Q_\thv$ is aperiodic; or
\item$g_\thv(x,y)>0$ for all $x,y$, and
$\nu$ has mass in each periodic
class of $Q_\thv$.
\end{enumerate}
Then ${\hat{\theta}}_{\nu,n} \stackrel{n\to\infty
}{\longrightarrow} [\thv]$,
$\PPs_\thv$-a.s.
\end{theorem}

The proof of this theorem is given in Section
\ref{sec:main-proof}.\vadjust{\goodbreak}
\begin{rem}
The assumptions 1--3 in Theorem \ref{thm:consistency} impose different
requirements on the initial measure $\nu$ used for the maximum
likelihood procedure. When the true parameter is aperiodic and has
nondegenerate observations, consistency holds for any choice of $\nu$.
On the other hand, in the case of degenerate observations, it is evident
that we cannot expect consistency to hold in general without imposing an
absolute continuity assumption of the form $\nu\sim\pi_\thv$. The
intermediate case, where the observations are nondegenerate but the signal
may be periodic, is not entirely obvious. An illuminating counterexample,
which shows that the MLE can be inconsistent for a choice of $\nu$ that
does not satisfy the requisite assumption in this case, is given in Remark
\ref{rem:periodicity} below.
\end{rem}
\begin{rem}
In Theorem \ref{thm:consistency}, we have assumed that the data is generated
by the stationary measure $\PPs_\thv$.
However, it follows directly from Lemma \ref{lem:ergodic} below
that, under the assumptions of Theorem \ref{thm:consistency}, we also have
${\hat{\theta}}_{\nu,n} \stackrel{n\to\infty}{\longrightarrow}
[\thv]$
$\PP_\thv^\rho\mbox{-a.s.}$
for any initial measure $\rho$ that satisfies the same assumptions
as $\nu$ in Theorem \ref{thm:consistency}. Hence, the initial measure
of the underlying chain is largely irrelevant, both for the consistency
of the estimator and in the definition of the log-likelihood function
$\ell_{\nu,n}(\theta)$ that is used to compute the estimator.
\end{rem}
\begin{rem}
The assumptions of Theorem \ref{thm:consistency} can be weakened
somewhat. For example, the $\sigma$-finite measure $\lambda$ can be
allowed to depend on $\theta$, or one may consider maximum likelihood
estimates of the form
$\hat\theta_n=\arg\max_{\theta\in\Theta}
\ell_{\nu_\theta,n}(\theta)$ where the initial measure $\nu$ used to
compute the likelihood depends on $\theta$ (the latter does not affect
the asymptotics of the MLE, but may improve finite sample properties in
certain cases). Such generalizations are straightforward and require
only minor adjustments in the proofs. In order not to further
complicate our notation, we leave these modifications to the reader.
\end{rem}
\begin{rem}
As was pointed out to us by a referee, assumptions
(A2) and (A4) depend on the choice of the observation reference
measure $\mu$, even though the maximum likelihood estimator itself is
independent of the choice of reference measure. It is therefore possible
that the assumptions of Theorem \ref{thm:consistency} are not satisfied
for a given reference measure $\mu$, but that consistency of the MLE can
be established nonetheless by making a suitable change of reference measure.
\end{rem}

\subsection{Geometric ergodicity implies identifiability}
\label{sec:sufficient-exponential-separation}

Most of the assumptions of Theorem \ref{thm:consistency} can be verified
in a straightforward manner. The exception is the identifiability
assumption (A6), which appears to be nontrivial.
Nonetheless, we will show that this assumption holds in a large class of
models: it is already sufficient (beside a mild technical assumption)
that the transition kernel $Q_\theta$ is \textit{geometrically
ergodic}, a
property that holds in many applications. Moreover, there is a
well-established theory of geometric ergodicity for Markov chains
\cite{meyntweedie1993} which provides a powerful set of tools to
verify this assumption. Consequently, our main theorem is directly
applicable in many cases of practical interest.
\begin{rem}
\label{rem:sketch-proof-v}
Before we state a precise result, it is illuminating to understand the
basic idea behind the proof of assumption (A6). Assume that
$Q_\theta$ is ergodic and that $\PPs_\theta^Y\ne\PPs_\thv^Y$. Then
there is an $s<\infty$ and a bounded function $h\dvtx\Yset^{s+1}\to\rset$
such that $\PEs_\theta[h(Y_0^s)]=0$ and $\PEs_\thv[h(Y_0^s)]=1$. Define
\[
A_n = \Biggl\{
y_1^n\dvtx \frac{1}{n-s}\sum_{i=1}^{n-s}
h(y_i^{i+s})>\frac{1}{2}
\Biggr\}
\]
for $n>s$. By the ergodic theorem, $\PPs_\thv^Y(A_n)\to1$ and
$\PPt^\lambda_\theta(A_n)\to0$ as $n\to\infty$. To prove
(A6), one must show that the convergence
$\PPt^\lambda_\theta(A_n)\to0$ happens at an exponential rate, that is,
one must establish a type of large deviations property. Therefore, the
key thing to prove is that geometrically ergodic Markov chains possess
such a large deviations property. This will be done in Section
\ref{sec:v-uniform-separation}.
\end{rem}

Let us begin by recalling the appropriate notion of geometric ergodicity
(the definition of the norm $\|\cdot\|_V$ was given in Section
\ref{sec:canonical} above).
\begin{defi}
\label{defi:v-uniformly-ergodic}
Let $V_\theta\dvtx\Xset\to[ 1,\infty)$ be given.
The transition kernel $Q_\theta$ is called
\textit{$V_\theta$-uniformly ergodic} if it possesses an invariant
probability measure $\pi_\theta$ and
\[
\|Q_\theta^m(x, \cdot)-\pi_\theta\|_{V_\theta}\le
R_\theta\alpha_\theta^{-m} V_\theta(x)\qquad
\mbox{for every }x\in\Xset, m\in\nset,
\]
for some constants $R_\theta<\infty$ and $\alpha_\theta>1$.
\end{defi}

For equivalent definitions and extensive discussion, see
\cite{meyntweedie1993}, Chapter 16. We can now formulate a
practical sufficient condition for assumption (A6).

\begin{enumerate}[(A6$'$)]
\item[(A6$'$)]
For every $\theta\not\sim\thv$ such that
$p^\lambda({Y}_{0}^{r_\theta};\theta)>0$ $\PPs_\thv\mbox{-a.s.}$,
there exists a function $V_\theta\ge1$ such
that $Q_\theta$ is $V_\theta$-uniformly ergodic,
$\PPs_\theta^Y\ne\PPs_\thv^Y$, and
%
%
\begin{equation}
\label{eq:conditional-lambda-integrable}
\PPs_\thv\biggl(
\int V_\theta(x_{r_\theta+1})
p^\lambda(\rmd x_{r_\theta+1},{Y}_{0}^{r_\theta};\theta)
<\infty
\biggr)>0 .
\end{equation}
\end{enumerate}
Note, in particular, that (\ref{eq:conditional-lambda-integrable}) holds
if (A4) holds and $|V_\theta|_\infty<\infty$ [in
this case, (A6$'$) implies that the transition kernel
$Q_\theta$ is \textit{uniformly ergodic}]. In the setting where
(A6$'$) holds, it is most natural to consider the equivalence
relation $\sim$ defined by setting $\theta\sim\theta'$ if and only if
$\PPs_\theta^Y=\PPs_{\theta'}^Y$ (i.e., two parameters are equivalent
precisely when they give rise to the same stationary observation laws).
\begin{theorem}
\label{thm:v-uniformly-ergodic}
Assume \textup{(A1)}, \textup{(A4)} and \textup{(A6$'$)}.
Then \textup{(A6)} holds.
\end{theorem}

The proof of this theorem is given in Section \ref{sec:proof-v-uniform}.

A different sufficient condition for assumption (A6), which
does not rely on geometric ergodicity of the underlying model, is the
following assumption (A6$''$). We will use this assumption in
Section \ref{sec:finite-state-leroux} to show that when $\Xset$ is
finite set, the identifiability assumption holds even for nonergodic
signals.
\begin{enumerate}[(A6$''$)]
\item[(A6$''$)] For every $\theta\not\sim\thv$ and
initial probability measure $\nu$, we have
\[
\liminf_{n\to\infty}\PPs_\thv(Y_0^n\in A_n)>0,\qquad
\limsup_{n\to\infty}n^{-1}\log\PP^\nu_\theta(Y_0^n\in A_n)<0
\]
for some sequence of sets $A_n\in\Ysigma^{\otimes(n+1)}$.
\end{enumerate}
\begin{prop}
\label{prop:mixing-for-all-initial}
Assume \textup{(A4)} and \textup{(A6$''$)}.
Then \textup{(A6)} holds.
\end{prop}

The proof of this proposition is given in Section \ref{sec:proof-v-uniform}.

\section{Examples}
\label{sec:applications}

In this section, we develop three classes of examples. In
Section~\ref{sec:postcomp:kalman} we consider linear Gaussian state space models.
In Section \ref{sec:finite-state-leroux}, we consider the classic case
where the signal state space is a finite set. Finally, in
Section~\ref{sec:example-nonlinear-state-space}, we develop a general class
of nonlinear state space models. In all these examples, we will find
that the assumptions of Theorem \ref{thm:consistency} are satisfied
in a rather general setting.

\subsection{Gaussian linear state space models}
\label{sec:postcomp:kalman}

Gaussian linear state space models form an important class of HMM. In
this setting, let $\Xset=\rset^d$ and $\Yset=\rset^p$ for some integers
$d,p$, and let $\Theta$ be a compact parameter space. The transition
kernel $T_\theta$ of the model is specified by the state space dynamics
%
%
\begin{eqnarray}
\label{eq:linear_state_space:time_depend:dynamic}
X_{k+1} & = & A_{\theta} X_k + R_{\theta} U_{k},
\\
\label{eq:linear_state_space:time_depend:observation}
Y_k & = &B_{\theta} X_k + S_{\theta} V_k,
\end{eqnarray}
where $\{(U_k,V_k)\}_{k \geq0}$ is an i.i.d. sequence of Gaussian
vectors with zero mean and identity covariance matrix, independent of $X_0$.
Here $U_k$ is $q$-dimensional, $V_k$ is $p$-dimensional, and the matrices
$A_\theta,R_\theta,B_\theta,S_\theta$ have the appropriate dimensions.

For each $\theta\in\Theta$ and any integer $r \geq1$, define
\[
\mathcal{O}_{\theta,r}\eqdef
\left[\matrix{ B_{\theta} \cr
B_{\theta} A_{\theta} \cr
B_{\theta} A_{\theta}^2 \cr
\vdots\cr
B_{\theta} A_{\theta}^{r-1}}\right] \quad\mbox{and}\quad
\mathcal{C}_{\theta,r}\eqdef
\left[\matrix{
R_{\theta} A_\theta R_{\theta} \cdots
A_\theta^{r-1}R_{\theta}}\right]
.
\]
It is assumed in the sequel that for any $\theta\in\Theta$, the following
hold:
\begin{longlist}[(L1)]
\item[(L1)] The pair $[A_\theta,B_\theta]$ is observable and the
pair $[A_\theta,R_{\theta}]$ is controllable, that is, the observability
matrix $\mathcal{O}_{\theta,d}$ and controllability matrix
$\mathcal{C}_{\theta,d}$ are full rank.
\item[(L2)]
The state transition matrix $A_\theta$ is discrete-time Hurwitz, that is,
its eigenvalues all lie in the open unit disc in $\mathbb{C}$.
\item[(L3)] The measurement noise covariance matrix $S_{\theta}$ is
full rank.
\item[(L4)] The functions $\theta\mapsto A_{\theta}$, $\theta
\mapsto
R_{\theta}$, $\theta\mapsto B_{\theta}$ and $\theta\mapsto
S_{\theta}$ are continuous on~$\Theta$.
\end{longlist}
We show below that the Markov kernel $Q_\theta$ is ergodic for every
$\theta\in\Theta$. We can therefore define without ambiguity the
equivalence relation $\sim$ on $\Theta$ as follows: $\theta\sim
\thetapr$
iff $\PPs_\theta^Y=\PPs_\thetapr^Y$.
We now proceed to verify the assumptions of Theorem \ref{thm:consistency}.

The fact that $A_\theta$ is Hurwitz guarantees that the state equation
is stable. Together with the controllability assumption, this implies
that $Q_\theta$ is $V_\theta$-uniformly ergodic with $V_\theta(x)
\asymp
|x|^2$ as $|x| \to\infty$ (\cite{glynnmeyn1996}, pages 929 and 930).
In particular, $Q_\thv$ is $V_\thv$-uniformly ergodic, which implies
(A1).

By the assumption that $S_\theta$ is full rank, and choosing the
reference measure $\mu$ to be the Lebesgue measure on $\Yset$,
we find that $g_\theta(x,y)$ is a Gaussian density for each $x\in
\Xset$
with covariance matrix $S_\theta S_\theta^T$. We therefore have
$|g_{\thv}|_\infty=(2\pi)^{-p/2}\det{}^{-1/2}(S_\thv S_\thv
^T)<\infty$,
so that $\PEs_\thv[\sup_{x\in\Xset} (\log g_\thv
(x,Y_0))^+]<\infty$.
On the other hand, as the stationary distribution $\pi_\theta$ is
Gaussian, the function $y\mapsto\int g_\thv(x,y) \pi_\thv(\rmd x)$
is a
Gaussian density with respect to $\mu$. Therefore, is easily seen that
$\PEs_\thv[|{\log\int g_\thv(x,Y_0) \pi_\thv(\rmd
x)}|]
<\infty$, and we have established (A2).

The dimension $q$ of the state noise vector $U_k$ is in many situations
smaller than the dimension $d$ of the state vector $X_k$ and hence
$R_{\theta}R_{\theta}^{T}$ may be rank deficient. However, note that
$Q_\theta^d(x,\rmd x')$ is a Gaussian distribution with covariance matrix
$\mathcal{C}_{\theta,d}\mathcal{C}_{\theta,d}^T$ for each $x\in
\Xset$.
Therefore, the controllability of the pair
$[A_\theta,R_\theta]$ nonetheless guarantees that $Q_\theta^d(x,\rmd x')$
has a density with respect to the Lebesgue measure $\lambda$ on $\Xset$.
Thus, (A3) is satisfied with $l=d$.

To proceed, we obtain an explicit expression for $p^\lambda
(y_0^r;\theta)$.
\begin{lem}
\label{lem:gaussian-p-lambda}
For $r\ge d$, we have
%
%
\begin{eqnarray}
\label{eq:expression-p-lambda}
p^\lambda({y}_{0}^{r-1};\theta)
&=&(2\pi)^{(d-pr)/2}
\det{}^{-1/2}(\mathcal{O}_{\theta,r}^T \Gamma_{\theta,r}^{-1}
\mathcal{O}_{\theta,r})
\det{}^{-1/2}(\Gamma_{\theta,r})\nonumber\\[-8pt]\\[-8pt]
&&{}\times
\exp\bigl( -\tfrac{1}{2} \mathbf{y}_r^T H_{\theta,r} \mathbf{y}_r
\bigr) .\nonumber
\end{eqnarray}
Here we defined the matrix
$\Gamma_{\theta,r}\stackrel{\mathit{def}}{=}\mathcal{H}_{\theta,r}\mathcal{H}_{\theta,r}^T
+\mathcal{S}_{\theta,r}\mathcal{S}_{\theta,r}^T$ with
\[
\mathcal{H}_{\theta,r} \stackrel{\mathit{def}}{=}
\pmatrix{
0 & 0 & \cdots& 0 \cr
B_\theta R_{\theta} & 0 & & 0 \cr
B_\theta A_\theta R_{\theta} & B_\theta R_{\theta}
& \cdots& 0 \cr
\vdots&\vdots&\ddots&\vdots\cr
B_\theta A_\theta^{r-2} R_{\theta} & B_\theta
A_\theta^{r-3} R_\theta& \cdots& B_\theta R_{\theta}}
\]
and where $\mathcal{S}_{\theta,r}$ is the $pr \times pr$ block diagonal
matrix with diagonal blocks equal to $S_{\theta}$,
$\mathbf{y}_r=[y_0,\ldots,y_{r-1}]^T$,
and $H_{\theta,r}$ is the matrix defined by
\[
H_{\theta,r}
\stackrel{\mathit{def}}{=}\Gamma_{\theta,r}^{-1} - \Gamma_{\theta,r}^{-1}
\mathcal{O}_{\theta,r}(\mathcal{O}_{\theta,r}^T \Gamma_{\theta,r}^{-1}
\mathcal{O}_{\theta,r})^{-1} \mathcal{O}_{\theta,r}^T\Gamma
_{\theta,r}^{-1}.
\]
\end{lem}
\begin{pf}
Define the vectors $\mathbf{Y}_r=[Y^T_0,\ldots, Y^T_{r-1}]^T$,
$\mathbf{U}_{r-1}= [U^T_0,\ldots,U^T_{r-2}]^T$ and
$\mathbf{V}_r=[V^T_0,\ldots, V^T_{r-1}]^T$.
It follows from elementary algebra that
\[
\mathbf{Y}_r=
\mathcal{O}_{\theta,r}X_0+
\mathcal{H}_{\theta,r}\mathbf{U}_{r-1}+
\mathcal{S}_{\theta,r}\mathbf{V}_r
\]
for any integer $r \geq1$. Note that, as $\mathbf{U}_{r-1}$ and
$\mathbf{V}_r$ are independent, the covariance matrix of the vector
$\mathcal{H}_{\theta,r}\mathbf{U}_{r-1}+\mathcal{S}_{\theta
,r}\mathbf{V}_r$
is given by $\Gamma_{\theta,r}$. It follows that
\[
p^x(y_0^{r-1};\theta) =
(2\pi)^{-pr/2} \det{}^{-1/2}(\Gamma_{\theta,r})
\exp\bigl( -\tfrac{1}{2}
(\mathbf{y}_r-\mathcal{O}_{\theta,r}x)^T \Gamma_{\theta,r}^{-1}
(\mathbf{y}_r-\mathcal{O}_{\theta,r}x)
\bigr) ,
\]
where we have used that $\Gamma_{\theta,r}$ is positive definite
(this follows directly from the assumption that $S_\theta$ is full rank).
Now let $\tilde{\Pi}_{\theta,r} \eqdef
\tilde{\mathcal{O}}_{\theta,r} ( \tilde{\mathcal
{O}}^T_{\theta,r}
\tilde{\mathcal{O}}_{\theta,r} )^{-1} \tilde{\mathcal
{O}}_{\theta,r}^T$
be the orthogonal projector on the range of $\tilde{\mathcal
{O}}_{\theta,r}
\eqdef\Gamma_{\theta,r}^{-1/2} \mathcal{O}_{\theta,r}$ ($\tilde
{\Pi
}_{\theta,r}$
is well defined for $r\ge d$ as the pair $[A_\theta,B_\theta]$
is observable, so that $\tilde{\mathcal{O}}_{\theta,r}$ is full rank).
Clearly,
\begin{eqnarray*}
(\mathbf{y}_r-\mathcal{O}_{\theta,r}x)^T \Gamma_{\theta,r}^{-1}
(\mathbf{y}_r-\mathcal{O}_{\theta,r}x) &=&
\|\tilde{\Pi}_{\theta,r}\Gamma_{\theta,r}^{-1/2}
\mathbf{y}_r - \Gamma_{\theta,r}^{-1/2}\mathcal{O}_{\theta,r}x
\|^2\\
&&{}+ \|(1-\tilde{\Pi}_{\theta,r})\Gamma_{\theta,r}^{-1/2}
\mathbf{y}_r\|^2.
\end{eqnarray*}
The result now follows from
\[
\int\exp\biggl(-\frac{1}{2}
\|\tilde{\Pi}_{\theta,r}\Gamma_{\theta,r}^{-1/2}
\mathbf{y}_r - \Gamma_{\theta,r}^{-1/2}\mathcal{O}_{\theta,r}x
\|^2
\biggr) \,\rmd x =
(2\pi)^{d/2}\det{}^{-1/2}(\mathcal{O}_{\theta,r}^T
\Gamma_{\theta,r}^{-1}\mathcal{O}_{\theta,r})
\]
(which is immediately seen to be finite due to the fact that
$\tilde{\mathcal{O}}_{\theta,r}$ has full rank), and from the
identity $H_{\theta,r}=\Gamma_{\theta,r}^{-1/2}
(1-\tilde{\Pi}_{\theta,r})\Gamma_{\theta,r}^{-1/2}$.
\end{pf}
\begin{rem}
As is evident from the proof, the observability assumption is key in
order to guarantee that $p^\lambda(y_0^{r-1};\theta)$ is finite
(albeit only
for $r$ sufficiently large). Intuitively, observability guarantees that
we can estimate $X_0$ from $Y_0^{d-1}$ ``in every direction,'' so that
the likelihood $p^x(y_0^{r-1};\theta)$ becomes small as $|x|\to\infty$.
This is needed in order to ensure that $p^x(y_0^{r-1};\theta)$ is
integrable with respect to the $\sigma$-finite measure $\lambda$. It
should also be noted that for any $r \geq d$ the matrix $H_{\theta,r}$
is rank-deficient, showing that (\ref{eq:expression-p-lambda}) is not
the density of a finite measure.
\end{rem}

Now note that, by our assumptions, the functions $\theta\mapsto
\det{}^{-1/2}(\mathcal{O}_{\theta,d}^T \Gamma_{\theta,d}^{-1}
\mathcal{O}_{\theta,d})$, $\theta\mapsto
\det{}^{-1/2}(\Gamma_{\theta,d})$, and $\theta\mapsto
H_{\theta,r}$ are continuous on $\Theta$ for any $r\ge d$.
Thus, $\theta\mapsto p^\lambda(y_0^{r-1};\theta)$ is continuous for
every $r\ge d$, and it is easily established that
$\PEs_\thv[\sup_{\thetapr\in\calU_\theta}
(\log p^\lambda({Y}_{0}^{r_\theta}; \thetapr))^+]
<\infty$ if we choose $r_\theta=d-1$ and a sufficiently small
neighborhood $\mathcal{U}_\theta$.
Moreover, note that
$|g_{\theta}|_\infty=(2\pi)^{-p/2}\det{}^{-1/2}(S_\theta
S_\theta
^T)$ and
$|q_{\theta}|_\infty=(2\pi)^{-d/2}\det{}^{-1/2}(\mathcal
{C}_{\theta,d}
\mathcal{C}_{\theta,d}^T)$. Therefore, by the continuity of
$S_\theta$ and $\mathcal{C}_{\theta,d}$, we have
$\sup_{\thetapr\in\calU_\theta}|q_{\thetapr}|_\infty<\infty$ and
$\PEs_\thv[\sup_{\thetapr\in\calU_\theta}
\sup_{x\in\Xset} (\log g_\thetapr(x,Y_0))^+]<\infty$
for a sufficiently small neighborhood $\mathcal{U}_\theta$.
Thus, we have verified (A4) and (A5).

It remains to establish assumption (A6). We established
above that $Q_\theta$ is $V_\theta$-uniformly ergodic with $V_\theta(x)
\asymp|x|^2$ as $|x| \to\infty$. Moreover, $\theta\not\sim\thv$ implies
$\PPs_\theta^Y\ne\PPs_\thv^Y$ by definition. Therefore, (A6$'$)
would be established if
\[
\int|x_d|^2 p^\lambda(\rmd x_d,Y_0^{d-1};\theta)<\infty,\qquad
\PPs_\thv\mbox{-a.s.}
\]
But note that
\[
\int p^\lambda(\rmd x_d,Y_0^{d-1};\theta) =
p^\lambda(Y_0^{d-1};\theta)<\infty,\qquad\PPs_\thv\mbox{-a.s.},
\]
so that $p^\lambda(\rmd x_d,Y_0^{d-1};\theta)$ is a finite measure. Moreover,
as $(Y_0^{d-1},X_d)=MX_0+\xi$ for a matrix $M$ and a Gaussian vector
$\xi$,
it is easily seen that $p^\lambda(\rmd x_d,Y_0^{d-1};\theta)$ must be
a random Gaussian measure. As Gaussian measures have finite moments,
we have established (A6$'$). Therefore,
(A6) follows from Theorem \ref{thm:v-uniformly-ergodic}.

Having verified (A1)--(A6), we can apply
Theorem \ref{thm:consistency}. As $g_\thv(x,y)>0$ for all $x,y$, and as
$Q_\thv$ is $V_\thv$-uniformly ergodic (hence certainly aperiodic), we
find that the MLE is consistent for any initial measure $\nu$.

\subsection{Finite state models}
\label{sec:finite-state-leroux}

One of the most widely used classes of HMM is obtained when the signal
is a finite state Markov chain. In this setting, let
$\Xset=\{1,\ldots,d\}$ for some integer $d$, let $\Yset$ be any Polish
space, and let $\Theta$ be a compact metric space. For each parameter
$\theta\in\Theta$, the signal transition kernel $Q_\theta$ is determined
by the corresponding transition probability matrix $\mathcal{Q}_\theta$,
while the observation density $g_\theta$ is given as in the general
setting of this paper.

It is assumed in the sequel that:
\begin{longlist}[(F1)]
\item[(F1)] The stochastic matrix $\mathcal{Q}_\thv$ is irreducible.
\item[(F2)] $\PEs_\thv[|{\log g_\thv(x,Y_0)}|]<\infty$ for every $x\in
\Xset$.
\item[(F3)] For every $\theta\in\Theta$, there is a neighborhood $\mathcal
{U}_\theta$
of $\theta$ such that
\[
\PEs_\thv\Bigl[\sup_{\theta'\in\mathcal{U}_\theta}
(\log g_{\theta'}(x,Y_0))^+\Bigr]<\infty\qquad
\mbox{for all }x\in\Xset.
\]
\item[(F4)] $\theta\mapsto\mathcal{Q}_\theta$ and $\theta\mapsto
g_\theta(x,y)$
are continuous for any $x\in\Xset$, $y\in\Yset$.
\end{longlist}
Following \cite{leroux1992}, we introduce the equivalence relation on
$\Theta$ as follows: we write $\theta\sim\theta'$ iff there exist
invariant distributions $\pi,\pi'$ for
$\mathcal{Q}_\theta,\mathcal{Q}_{\theta'}$, respectively, such that
$\PP_\theta^{\pi,Y}=\PP_{\theta'}^{\pi',Y}$. In words, two parameters
are equivalent whenever they give rise to the same stationary
observation laws for some choice of invariant measures for the underlying
signal process. The latter statement is not vacuous as we have not
required that $\mathcal{Q}_\theta$ is ergodic for $\theta\ne\thv$, that
is, there may be multiple invariant measures for $\mathcal{Q}_\theta$.
The possibility that $\mathcal{Q}_\theta$ is not aperiodic or even
ergodic is the chief complication in this example, as the easily
verified $V$-uniform ergodicity assumption (A6$'$) need not
hold. We will show nonetheless that assumption (A6$''$) is
satisfied, so that Theorem \ref{thm:consistency} can be applied.
\begin{lem}
\label{lem:finite-state-hoeffding}
Let $C\subseteq\Xset$ be an ergodic class of $\mathcal{Q}_\theta$, and
denote by $\pi_C$ the unique $\mathcal{Q}_\theta$-invariant measure
supported in $C$.
Fix $s \geq0$, and let $f\dvtx \Yset^{s+1} \to\rset$ be such that
$|f|_\infty<\infty$. Then there exists a constant $K$ such that
\[
\PP_\theta^\nu\Biggl( \Biggl|\sum_{i=1}^n
\{f({Y}_{i}^{i+s}) - \PP_\theta^{\pi_C}[f({Y}_{0}^{s})]
\}
\Biggr| \geq t \Biggr) \leq
K \exp\biggl[ - \frac{t^2}{Kn} \biggr]
\]
for any probability measure $\nu$ supported in $C$ and any $t > 0$,
$n\ge1$.
\end{lem}
\begin{pf}
The proof is identical to that of Theorem
\ref{thm:azuma-hoeffding-v-unif}, provided we replace the application of
Theorem \ref{thm:exponential-inequality} by a trivial modification of
the result of \cite{glynnormoneit2002}.
\end{pf}
\begin{rem}
As stated, the result of \cite{glynnormoneit2002} would require that
the restriction of $\mathcal{Q}_\theta$ to $C$ is aperiodic. However,
aperiodicity is only used in the proof to ensure the existence of a
solution to the Poisson equation, and it is well known that the latter
holds also in the periodic case. Therefore, a trivial modification of
the proof in \cite{glynnormoneit2002} allows us to apply the result
without additional assumptions.
\end{rem}
\begin{lem}
\label{lem:finite-state-mixing-check}
In the present setting, assumption \textup{(A6$''$)} holds.
\end{lem}
\begin{pf}
Let $\theta\not\sim\thv$. We can partition
$\Xset=E_1\cup\cdots\cup E_p\cup T$ into the $p\le d$ ergodic classes
$E_1,\ldots,E_p$ and the set of transient states $T$ of the stochastic
matrix~$\mathcal{Q}_\theta$. Denote as $\pi_\theta^i$ the unique invariant
measure of $\mathcal{Q}_\theta$ that is supported in~$E_i$. Then we can
find an integer $s\ge1$ and bounded function
$h\dvtx\Yset^{s+1}\to\mathbb{R}$ such that $\PE_\theta^{\pi_\theta^i}[
h(Y_0^s)]\le0$ for all $i=1,\ldots,p$ and such that
$\PEs_\thv[h(Y_0^s)]=1$.

Define for $n>2s$ the set $A_n\in\Ysigma^{\otimes(n+1)}$ as
\[
A_n \eqdef
\Biggl\{{y}_{0}^{n} \in\Yset^{n+1}\dvtx
\frac{1}{\lfloor n/2\rfloor-s}\sum_{i=\lceil n/2\rceil+1}^{n-s}
h({y}_{i}^{i+s}) \geq\frac{1}{2} \Biggr\}.
\]
As $Y_0^\infty$ is stationary and ergodic under $\PPs_\thv$ (because
$\mathcal{Q}_\thv$ is irreducible), we have
\[
\lim_{n\to\infty}\PPs_\thv(Y_0^n\in A_n) =
\lim_{n\to\infty}\PPs_\thv \Biggl[
\frac{1}{\lfloor n/2\rfloor-s}\sum_{i=1}^{\lfloor n/2\rfloor-s}
h({Y}_{i}^{i+s}) \geq\frac{1}{2}
\Biggr] = 1
\]
by Birkhoff's ergodic theorem. On the other hand, for any initial
probability measure $\nu$, we can estimate as follows: for some constant
$K>0$,
\begin{eqnarray*}
\PP_\theta^\nu(Y_0^n\in A_n)
&=&
\PP_\theta^\nu\bigl(Y_0^n\in A_n, X_{\lceil n/2\rceil}\in T\bigr)\\
&&{} +
\sum_{j=1}^p
\PP_\theta^\nu\bigl(Y_0^n\in A_n|X_{\lceil n/2\rceil}\in E_j\bigr)
\PP_\theta^\nu\bigl(X_{\lceil n/2\rceil}\in E_j\bigr)
\\
&\le&
\PP_\theta^\nu\bigl(X_{\lceil n/2\rceil}\in T\bigr)\\
&&{} +
\max_{j=1,\ldots,p}
\sup_{\operatorname{supp}\mu\subseteq E_j}\PP_\theta^\mu
\Biggl[\frac{1}{\lfloor n/2\rfloor-s}\sum_{i=1}^{\lfloor n/2\rfloor-s}
h({Y}_{i}^{i+s}) \geq\frac{1}{2} \Biggr]
\\
&\le&
K \exp\biggl[-\frac{n}{K}\biggr].
\end{eqnarray*}
The latter inequality follows from the fact that the population in the
transient states decays exponentially, while we may apply Lemma
\ref{lem:finite-state-hoeffding} to obtain an exponential bound
for every ergodic class $E_i$. We therefore find that
\[
\limsup_{n\to\infty}n^{-1}\log\PP^\nu_\theta(Y_0^n\in A_n)\le
-\frac{1}{K}<0,
\]
completing the proof of assumption (A6$''$).
\end{pf}

Let us now check the assumptions of Theorem
\ref{thm:consistency}. (A1) follows
directly from the assumption that $\mathcal{Q}_\thv$ is irreducible.
To establish (A2), note that
\[
\PEs_\thv\Bigl[\sup_{x\in\Xset} (\log g_\thv(x,Y_0))^+\Bigr]
\le
\sum_{x\in\Xset}\PEs_\thv[|{\log g_\thv(x,Y_0)}|]
< \infty,
\]
while we can estimate
\begin{eqnarray*}
&&\PEs_\thv\biggl[\biggl|{\log\int g_\thv(x,Y_0) \pi_\thv(\rmd x)}
\biggr|\biggr]\\
&&\qquad\le
\PEs_\thv\Bigl[\sup_{x\in\Xset} (\log g_\thv(x,Y_0))^+\Bigr]
+ \PEs_\thv\Bigl[\sup_{x\in\Xset} (\log g_\thv(x,Y_0))^-\Bigr]
\\
&&\qquad\le
\sum_{x\in\Xset}\PEs_\thv[|{\log g_\thv(x,Y_0)}|] <
\infty.
\end{eqnarray*}
Assumption (A3) holds trivially for $l=1$ and with $\lambda$ the
counting measure on $\Xset$ [note that $|q_{\theta}|_\infty\le1$
for all
$\theta$, as $q_\theta(x,x')$ is simply the transition probability
from $x$ to $x'$]. To establish (A4), note that
$\sup_{\theta\in\Theta}|q_{\theta}|_\infty<\infty$, while
\[
\PEs_\thv\Bigl[\sup_{\thetapr\in\calU_\theta}
\sup_{x\in\Xset} (\log g_\thetapr(x,Y_0))^+\Bigr]
\le
\sum_{x\in\Xset}
\PEs_\thv\Bigl[\sup_{\thetapr\in\calU_\theta}
(\log g_\thetapr(x,Y_0))^+\Bigr]
<\infty
\]
by our assumptions. Moreover, as
\[
\PEs_\thv\Bigl[\sup_{\thetapr\in\calU_\theta}
(\log p^\lambda({Y}_{0}^{0}; \thetapr))^+\Bigr]
\le
\PEs_\thv\Bigl[\sup_{\thetapr\in\calU_\theta}
\sup_{x\in\Xset}
\bigl(\log d + \log g_\thetapr(x,Y_0)\bigr)^+\Bigr]
<\infty,
\]
we have shown that (A4) holds with $r_\theta=0$ for
all $\theta$. Next, we note that the continuity of
$\theta\mapsto\mathcal{Q}_\theta$ and $\theta\mapsto g_\theta(x,y)$
yield immediately that $\theta\mapsto p^\lambda(y_0^n;\theta)$ is
a continuous function for every $n\ge0$ and $y_0^n\in\Yset^{n+1}$,
establishing (A5). Finally, Lemma
\ref{lem:finite-state-mixing-check} and Proposition
\ref{prop:mixing-for-all-initial} establish (A6).

Having verified (A1)--(A6), we can apply
Theorem \ref{thm:consistency}. Note that as $\mathcal{Q}_\thv$ is
irreducible, $\pi_\thv$ charges every point of $\Xset$. Therefore,
by Theorem \ref{thm:consistency}, the MLE is consistent provided that
$\nu$ charges every point of $\Xset$ (so that $\nu\sim\pi_\thv$).
\begin{rem}
The result obtained in this section as a special case of Theorem
\ref{thm:consistency} is almost identical to the result of Leroux
\cite{leroux1992}. The main difference in \cite{leroux1992} is that
there the parameter space $\Theta$ may be noncompact, provided the
parametrization of the model vanishes at infinity. This setting reduces
directly to the compact case by compactifying the parameter space
$\Theta$,
so that this does not constitute a major generalization from the
technical point of view.

However, it should be noted that one cannot immediately apply Theorem
\ref{thm:consistency} to the compactified model. The problem is that
the new parameters ``at infinity'' are typically sub-probabilities
rather than true probability measures, while we have assumed in this
paper that every parameter $\theta\in\Theta$ corresponds to a
probability measure on the space of observation paths. Theorem
\ref{thm:consistency} can certainly be generalized to allow for
sub-probabilities without significant technical complications.
We have chosen to concentrate on the compact setting, however, in order
to keep the notation and results of the paper as clean as possible.
\end{rem}

\subsection{Nonlinear state space models}
\label{sec:example-nonlinear-state-space}

In this section, we consider a class of nonlinear state space models.
Let $\Xset=\rset^d$, $\Yset=\rset^\ell$, and let $\Theta$ be a compact
metric space. For each $\theta\in\Theta$, the Markov kernel
$Q_\theta$
of the hidden process $(X_k)_{k\ge0}$ is defined through the nonlinear
recursion
\[
X_{k} = G_\theta(X_{k-1})+ \Sigma_\theta(X_{k-1}) \zeta_k .
\]
Here $(\zeta_k)_{k \geq1}$ is an i.i.d. sequence of $d$-dimensional
random vectors which are assumed to possess a density $\rho_\zeta$ with
respect to the Lebesgue measure $\lambda$ on $\rset^d$, and
$G_\theta\dvtx
\rset^d \to\rset^d$, $\Sigma_\theta\dvtx \rset^d \to\rset^{d \times d}$
are given (measurable) functions. The model for the hidden chain
$(X_k)_{k \geq0}$ is sometimes known as a vector ARCH model, and covers
many models of interest in time series analysis and financial
econometrics (including the AR model, the ARCH model, threshold ARCH,
etc.). We let the reference measure $\mu$ be the Lebesgue measure on
$\rset^\ell$, and define the observed process $(Y_k)_{k\ge0}$ by means
of a given observation density $g_\theta(x,y)$.

For any positive matrix $B$, denote by $\lambda_{\min}(B)$ its minimal
eigenvalue. For any bounded set $\mathcal{A} \subset
\rset^{d \times d}$, define $\mathcal{A}_m \eqdef\{ A_1A_2
\cdots
A_m\dvtx A_i \in\mathcal{A}, i=1,\ldots,m\}$. Denote by
$\rho(\mathcal{A})$ the joint spectral radius of the set of matrices
$\mathcal{A}$, defined as
\[
\rho(\mathcal{A} ) \eqdef
\limsup_{m \to\infty} \Bigl( {\sup_{A \in\mathcal{A}_m}} \| A \|
\Bigr)^{1/m} .
\]
Here \mbox{$\| \cdot\|$} is any matrix norm [it is elementary that
$\rho(\mathcal{A})$ does not depend on the choice of the norm].
We now introduce the basic assumptions of this section.
\begin{longlist}[(NL1)]
\item[(NL1)] 
The random variables $\zeta_k$ have mean zero and identity covariance
matrix. Moreover, $\rho_\zeta(x)>0$ for all $x\in\rset^d$, and
$|\rho_\zeta|_{\infty} < \infty$.
\item[(NL2)] 
For each $\theta\in\Theta$, the function $\Sigma_\theta$ is
bounded on
compact sets, $\Sigma_\theta(x)=o(|x|)$ as $|x|\to\infty$, and
$0< \inf_{\thetapr\in\mathcal{U}_\theta}
\inf_{x \in\rset^d} \lambda_{\min}[\Sigma_\thetapr(x)
\Sigma^T_\thetapr(x) ] $
for a sufficiently small neighborhood $\mathcal{U}_\theta$ of $\theta$.
\item[(NL3)] 
For each $\theta\in\Theta$, the drift function $G_\theta$ has the form
\[
G_\theta(x)= A_\theta(x) x + h_\theta(x)
\]
for some measurable functions $A_\theta\dvtx\rset^d\to\rset^{d\times
d}$ and
$h_\theta\dvtx\rset^d\to\rset^d$. Moreover, we assume that
$G_\theta$ is bounded on compact sets, $h_\theta(x)= o(|x|)$ as
$|x|\to\infty$,
and that there exists $R_\theta> 0$ such that the set of matrices
$\mathcal{A}_{\theta} \eqdef\{ A_\theta(x) \dvtx x \in\rset^d, |x|
\geq
R_\theta\}$ is bounded and $\rho(\mathcal{A}_{\theta}) < 1$.
\item[(NL4)] 
For each $\theta\in\Theta$, there is a neighborhood $\mathcal
{U}_\theta$
of $\theta$ such that
\begin{eqnarray*}
\PEs_\thv\Bigl[\sup_{\thetapr\in\calU_\theta}
\sup_{x\in\Xset} (\log g_\thetapr(x,Y_0))^+\Bigr]&<&\infty,
\\
\PEs_\thv\biggl[\sup_{\thetapr\in\calU_\theta}
\biggl(\log\int g_\thetapr(x,Y_0) \lambda(\rmd x)\biggr)^+\biggr]
&<&\infty.
\end{eqnarray*}
Moreover, $\PPs_\thv(\int|x| g_\theta(x,Y_0) \lambda(\rmd
x)<\infty
)>0$ for each $\theta\in\Theta$, and
\[
\PEs_{\thv}\biggl[ \int(\log g_\thv(x,Y_0))^- \pi_{\thv}(\rmd x)
\biggr] < \infty.
\]
\item[(NL5)] 
The functions $\theta\mapsto g_\theta(x,y)$, $\theta\mapsto G_\theta(x)$,
$\theta\mapsto\Sigma_\theta(x)$ and $x\mapsto\rho_\zeta(x)$ are
continuous
on $\Theta$ for every $x,y$. Moreover, for each $\theta\in\Theta$,
the function $\thetapr\mapsto\int g_\thetapr(x,Y_0) \lambda(\rmd x)$
is positive and continuous at $\theta$, $\PPs_{\thv}$-a.s.
\end{longlist}
\begin{rem}
We have made no attempt at generality here: for the sake of example, we have
chosen a set of conditions under which the assumptions of
Theorem~\ref{thm:consistency} are easily verified. Of course, the applicability
of Theorem \ref{thm:consistency} extends far beyond the simple
assumptions imposed in this section.

Nonetheless, even the present assumptions already cover a broad class of
nonlinear models. Consider, for example, the stochastic volatility
model \cite{hullwhite1987}
%
%
\begin{equation}
\cases{
X_{k+1}= \phi_\theta X_k + \sigma_\theta \zeta_k, \cr
Y_k= \beta_\theta\exp( X_k/2) \varepsilon_k,}
\end{equation}
where $(\zeta_k,\varepsilon_k)$ are i.i.d. Gaussian random variables in
$\rset^2$ with zero mean and identity covariance matrix, $\beta
_\theta>
0$, $\sigma_\theta> 0$ and $|\phi_\theta| < 1$ for every
$\theta\in\Theta$, and the functions $\theta\mapsto\phi_\theta$,
$\theta\mapsto\sigma_\theta$ and $\theta\mapsto\beta_\theta$ are
continuous. Assumptions
(NL1)--(NL3)
are readily seen to hold. The observation likelihood
$g_\theta$ is given by
\[
g_\theta(x,y) = (2 \pi\beta_\theta^2)^{-1/2}
\exp[ - \exp(-x) y^2 /2\beta_\theta^2 - x/2] .
\]
We can compute
\[
\sup_{x\in\Xset}g_\theta(x,y) = \frac{1}{\sqrt{2\pi\rme}}
\frac{1}{|y|},\qquad
\int g_\theta(x,y) \lambda(\rmd x) = \frac{1}{|y|}.
\]
As the stationary distribution $\pi_{\thv}$ is Gaussian, it is easily
seen that the law of $Y_0$ under $\PPs_\thv$ has a bounded density
with respect to the Lebesgue measure $\mu$ on $\Yset$.
As $\int(\log(1/|y|))^+\mu(\rmd y)<\infty$, the first
equation\vspace*{1pt}
display of (NL4) follows. To prove that
$\PPs_\thv(\int|x| g_\theta(x,Y_0) \lambda(\rmd x)<\infty
)>0$, it suffices to note that $x\mapsto g_\theta(x,y)$ has
exponentially decaying tails for all $|y|>0$. The remaining part of
(NL4) follows easily using that
$\pi_\thv$ is Gaussian and $\PEs_\thv(Y_0^2)<\infty$. Finally,
(NL5) now follows immediately, and we have verified
that the assumptions of this section hold for the stochastic volatility
model. Similar considerations apply in a variety of nonlinear models
commonly used in financial econometrics.
\end{rem}

We show below that the Markov kernel $Q_\theta$ is ergodic for every
$\theta\in\Theta$. We can therefore define without ambiguity the
equivalence relation $\sim$ on $\Theta$ as follows: $\theta\sim
\thetapr$
iff $\PPs_\theta^Y=\PPs_\thetapr^Y$.
We now proceed to verify the assumptions of Theorem \ref{thm:consistency}.

It is shown in \cite{liebscher2005}, Theorem 2, that under conditions
(NL1)--(NL3), the Markov
kernel $Q_\theta$ is $V$-uniformly ergodic for each $\theta\in\Theta$
with $V(x)=1+|x|$. In particular, assumption
(A1)
holds. The first part of (A2) follows directly
from (NL4). To prove the second part,
we first note that $Q_\theta$ has a transition density
\[
q_\theta(x,x') =
|{\det }[\Sigma_\theta(x) ] |^{-1}
\rho_\zeta\bigl( \Sigma_\theta^{-1}(x) \{x'- G_\theta(x)\}
\bigr)
\]
with respect to the Lebesgue measure $\lambda$ on $\Xset$. This
evidently gives
\[
|q_{\theta}|_\infty=\sup_{x\in\Xset}|{\det}[\Sigma_\theta
(x)]|^{-1}
|\rho_\zeta|_\infty<\infty
\]
by (NL1) and (NL2),
which implies in particular that
$\pi_\thv$ has a bounded density with respect to $\lambda$. Therefore
\begin{eqnarray*}
&&\PEs_\thv\biggl[\biggl(\log\int g_\thv(x,Y_0) \pi_\thv(\rmd x)
\biggr)^+\biggr] \\
&&\qquad\le
|q_{\thv}|_\infty
\PEs_\thv\biggl[\biggl(\log\int g_\thv(x,Y_0) \lambda(\rmd x)
\biggr)^+\biggr] < \infty
\end{eqnarray*}
by (NL4). On the other hand,
as $x\mapsto(\log x)^-$ is convex, we have
\[
\PEs_\thv\biggl[\biggl(\log\int g_\thv(x,Y_0) \pi_\thv(\rmd x)
\biggr)^-\biggr] \le
\PEs_\thv\biggl[\int(\log g_\thv(x,Y_0))^- \pi_\thv(\rmd x)
\biggr] < \infty
\]
by Jensen's inequality and (NL4). Therefore,
(A2) is established. We have already shown that
$Q_\theta$ possesses a bounded density, so (A3) holds
with $l=1$. Assumption (A4) with $r_\theta=0$ follows
directly from (NL4) and
(NL1), (NL2).

To establish (A5), let $\nu_\theta(\rmd x,y)
\eqdef g_\theta(x,y) \lambda(\rmd x)/\int g_\theta(x,y) \lambda
(\rmd
x)$. By (NL5), $\nu_\theta(\rmd x,Y_0)$ is a
probability measure $\PPs_\thv$-a.s., and for every $\theta\in
\Theta$
the density function $\thetapr\mapsto g_\thetapr(x,Y_0)/\int
g_\thetapr(x,Y_0) \lambda(\rmd x)$ is continuous at $\theta$
$\PPs_\thv$-a.s. By Scheff\'e's lemma, this implies that for any
$\theta\in\Theta$, the map $\thetapr\mapsto\nu_\thetapr(\cdot
,Y_0)$ is
continuous at $\theta$ $\PPs_\thv$-a.s. with respect to the total
variation norm \mbox{$\|\cdot\|_1$}. Similarly, as $\theta\mapsto q_\theta(x,x')$
is continuous by (NL5), the map
$\theta\mapsto Q_\theta(x,\rmd x')$ is continuous with respect to the
total variation norm. Now note that we can write
\[
p^\lambda(Y_0^n;\theta) =
\biggl(
\int g_\theta(x,Y_0) \lambda(\rmd x)\biggr)
\int p^{x'}(Y_1^n;\theta) Q_\theta(x,\rmd x')
\nu_\theta(\rmd x,Y_0).
\]
From (NL4), it follows that $x\mapsto
\sup_{\thetapr\in\mathcal{U}_\theta}g_\thetapr(x,Y_k)$ is bounded
$\PPs_\thv$-a.s. for every~$k$. Therefore, $x\mapsto
\sup_{\thetapr\in\mathcal{U}_\theta}p^x(Y_1^n;\thetapr)$ is a bounded
function $\PPs_\thv$-a.s., and by dominated convergence the function
$\thetapr\mapsto p^x(Y_1^n;\thetapr)$ is continuous at $\theta$
$\PPs_\thv$-a.s. for every $\theta\in\Theta$. Therefore, it follows
that $\PPs_\thv$-a.s.
\begin{eqnarray*}
&&
\biggl|\int p^{x'}(Y_1^n;\theta_n) Q_{\theta_n}(x,\rmd x')
\nu_{\theta_n}(\rmd x,Y_0) -
\int p^{x'}(Y_1^n;\theta) Q_\theta(x,\rmd x')
\nu_\theta(\rmd x,Y_0)\biggr| \\
&&\qquad\le
\int|p^{x'}(Y_1^n;\theta_n)-p^{x'}(Y_1^n;\theta)|
Q_\theta(x,\rmd x') \nu_\theta(\rmd x,Y_0) \\
&&\qquad\quad{} +
\sup_{\thetapr\in\mathcal{U}_\theta}|p^{\cdot}(Y_1^n;\thetapr
)|_\infty
\|\nu_{\theta_n}(\cdot,Y_0)Q_{\theta_n}
-\nu_\theta(\cdot,Y_0)Q_\theta\|_1\\
&&\qquad
\stackrel{n\to\infty}{\longrightarrow}0
\end{eqnarray*}
for any sequence $(\theta_n)_{n\ge0}\subset\mathcal{U}_\theta$,
$\theta_n\to\theta$. Here we have used the dominated convergence
theorem to conclude convergence of the first term, and the continuity in
total variation established above for the second term.
(A5) follows.

It remains to establish assumption (A6). We established
above that $Q_\theta$ is $V$-uniformly ergodic with $V(x)=1+|x|$.
Moreover, $\theta\not\sim\thv$ implies $\PPs_\theta^Y\ne\PPs
_\thv^Y$ by
definition. Therefore, (A6$'$) would be established if
\[
\PPs_\thv\biggl(
\int|x'| Q_\theta(x,\rmd x')
g_\theta(x,Y_0) \lambda(\rmd x)<\infty\biggr)>0.
\]
But as $Q_\theta$ is $V$-uniformly ergodic, it follows that $Q_\theta
V\le\alpha_\theta V+K_\theta$ for some positive constants
$\alpha_\theta,K_\theta$ (\cite{meyntweedie1993}, Theorem 16.0.1).
Assumption (A6$'$) therefore follows from
(NL4), and
(A6) follows from Theorem \ref{thm:v-uniformly-ergodic}.

Having verified (A1)--(A6), we can apply
Theorem \ref{thm:consistency}. As $g_\thv(x,y)>0$ for all $x,y$, and as
$Q_\thv$ is $V$-uniformly ergodic (hence certainly aperiodic), we
find that the MLE is consistent for any initial measure $\nu$.
\begin{rem}
The assumption in (NL4) that
\[
\PEs_{\thv}\biggl[ \int(\log g_\thv(x,Y_0))^- \pi_{\thv}(\rmd x)
\biggr] < \infty
\]
is used to verify the second part of (A2). This
condition can be replaced by the following assumption:
there exists a set $D \in\Xsigma$ such that:
\begin{longlist}
\item$\PEs_\thv[(\log\int_D g_\thv(x,Y_0)
\lambda(\rmd x))^-]<\infty$, and
\item$\inf_{x,x' \in D} q_\thv(x,x')>0$.
\end{longlist}
The latter condition is sometimes easier to check.

To see that the result still holds under this modified condition,
note that
\[
\PEs_\thv\biggl[\biggl(\log\int g_\thv(x,Y_0) \pi_\thv(\rmd x)
\biggr)^+\biggr]
<\infty
\]
follows as above. On the other hand,
\begin{eqnarray*}
\int g_\thv(x,Y_0) \pi_\thv(\rmd x) &\ge&
\int_{D\times D} g_\thv(x',Y_0) q_\thv(x,x') \pi_\thv(\rmd x)
\lambda(\rmd x') \\
&\ge&\pi_\thv(D)
\inf_{x,x'\in D}q_\thv(x,x')
\int_{D} g_\thv(x,Y_0)\lambda(\rmd x).
\end{eqnarray*}
It follows from (i) that $\lambda(D)>0$, so that $\pi_\thv(D)>0$ also
(as $Q_\thv$, and therefore~$\pi_\thv$, has a positive density with
respect to $\lambda$). It now follows directly that also
$\PEs_\thv[(\log\int g_\thv(x,Y_0) \pi_\thv(\rmd
x))^-]
<\infty$, and the claim is established.
\end{rem}

\section{\texorpdfstring{Proof of Theorem
\protect\ref{thm:consistency}}{Proof of Theorem 1}}
\label{sec:main-proof}

The proof of Theorem \ref{thm:consistency} consists of three parts.
First, we prove pointwise convergence of the log-likelihood under the
true parameter $\thv$ (Section \ref{sec:pointwise}). Next, we establish
identifiability of every $\theta\not\sim\thv$ (Section
\ref{sec:identifiability}). Finally, we put everything together to
complete the proof of consistency (Section~\ref{sec:consistency}).

\subsection{Pointwise convergence of the normalized log-likelihood under
$\thv$}
\label{sec:pointwise}

The goal of this section is to show that our hidden Markov model
possesses a finite entropy rate and that the asymptotic equipartition
property holds. We begin with a simple result, which will be used
to reduce to the stationary case.
\begin{lem}
\label{lem:ergodic}
Assume \textup{(A1)}. Then
$(Y_k)_{k\ge0}$ is ergodic under $\PPs_\thv$. Moreover, if one
of the assumptions 1--3 of Theorem \ref{thm:consistency} hold,
then $\PP_\thv^{\nu,Y}\sim\PPs_\thv^Y$.
\end{lem}
\begin{pf}
As $Q_\thv$ possesses a unique invariant measure by
(A1), the kernel $T_\thv$ possesses a unique
invariant measure also. This implies that the process $(X_k,Y_k)_{k\ge
0}$ is ergodic under the stationary measure $\PPs_\thv$ (as the latter
is trivially an extreme point of the set of stationary measures).
Therefore, $(Y_k)_{k\ge0}$ is ergodic also.

If $\nu\sim\pi_\thv$, it is easily seen that
$\PP_\thv^{\nu,Y}\sim\PPs_\thv^Y$ [as $d\PP_\thv^{\nu}/d\PPs
_\thv=
(d\nu/d\pi_\thv)(X_0)$]. Otherwise, we argue as follows.
Suppose that $Q_\thv$ has period $d$ [this is guaranteed to hold for some
$d\in\nset$ by (A1)]. Then we can partition the
signal state space as $\Xset= C_1\cup\cdots\cup C_d\cup F$, where
$C_1,\ldots,C_d$ are the periodic classes and $\pi(F)=0$
(\cite{meyntweedie1993}, Section 5.4.3). Note that
$C_1,\ldots,C_d$ are absorbing sets for $Q_\thv^d$ where the restriction
of $Q_\thv^d$ to $C_i$ is positive Harris and aperiodic with the
corresponding invariant probability measure $\pi_\thv^i$. Moreover, the
Harris recurrence assumption guarantees that
$\PP_\thv^x(X_n\notin F\mbox{ eventually})=1$ for all $x\in\Xset$.
Therefore, $\nu Q_\thv^{nd}(F)\downarrow0$ and
$\nu Q_\thv^{nd}(C_i)\uparrow c_\nu^i$ as $n\to\infty$. It follows
from the ergodic theorem for aperiodic Harris recurrent Markov chains that
\[
\|\nu Q_\thv^n - \pi_\thv^\nu Q_\thv^n\|_1
\stackrel{n\to\infty}{\longrightarrow},\qquad
\pi_\thv^\nu\eqdef\sum_{i=1}^dc_\nu^i\pi_\thv^i.
\]
Using $g_\thv(x,y)>0$ and \cite{vanhandel2009}, Lemma 3.7,
this implies that $\PP_\thv^{\nu,Y}\sim\PP_\thv^{\pi_\thv^\nu,Y}$.
But if $\nu$ has mass in each periodic class $C_i$ or if $d=1$, then
$c_\nu^i>0$ for all $i=1,\ldots,d$. Thus, $\pi_\thv^\nu\sim\pi=
\frac{1}{d}\sum_{i=1}^d\pi_\thv^i$ which
implies $\PPs_\thv^Y\sim\PP_\thv^{\pi_\thv^\nu,Y}\sim\PP_\thv
^{\nu,Y}$.
\end{pf}

We will also need the following lemma.
\begin{lem}
\label{lem:finite-entropy-for-all-n}
Assume \textup{(A2)}. Then
$\PEs_\thv[|{\log\pstat(Y_0^n;\thv)}|]<\infty$ for all $n\ge0$.
\end{lem}
\begin{pf}
We easily obtain the upper bound
\[
\PEs_\thv[(\log\pstat(Y_0^n;\thv))^+] \le
\PEs_\thv\Biggl[
\sum_{k=0}^n\sup_{x\in\Xset} (\log g_\thv(x,Y_k))^+
\Biggr] < \infty.
\]
On the other hand, we can estimate
\begin{eqnarray*}
\PEs_\thv[\log\pstat(Y_0^n;\thv)] &=&
\PEs_\thv\biggl[\log\frac{\pstat(Y_0^n;\thv)}{
\prod_{k=0}^n\int g_\thv(x,Y_k) \pi_\thv(\rmd x)
}\biggr] \\
&&{}+
\PEs_\thv\Biggl[\sum_{k=0}^n\log
\int g_\thv(x,Y_k) \pi_\thv(\rmd x)\Biggr] \\
&\ge&
-(n+1)
\PEs_\thv\biggl[\biggl|{\log\int g_\thv(x,Y_0) \pi_\thv(\rmd x)}
\biggr|\biggr]\\
&>& -\infty,
\end{eqnarray*}
where we have used that relative entropy is nonnegative.
\end{pf}

The main result of this section follows. After a reduction to the
stationary case by means of the previous lemma, the proof concludes by
verifying the assumptions of the generalized Shannon--Breiman--McMillan
theorem \cite{barron1985}.
\begin{theorem}
\label{theo:ShannonBreimanMcMillan}
Assume \textup{(A1)} and \textup{(A2)}.
There exists $-\infty<\ell(\thv)<\infty$ such that
%
%
\begin{equation}
\label{eq:mean-convergence-likelihood-mu}
\ell(\thv)= \lim_{n \to\infty} \PEs_\thv[ n^{-1} \log
\pstat(Y_0^n;\thv) ] ,
\end{equation}
and such that
%
%
\begin{equation}
\label{eq:as-convergence-likelihood-mu}
\ell(\thv)= \lim_{n \to\infty} n^{-1} \log p^\nu(Y_0^n;\thv
), \qquad\PPs_{\thv}\mbox{-a.s.}
\end{equation}
for any probability measure $\nu$ such that one
of the assumptions 1--3 of Theorem \ref{thm:consistency} is
satisfied (in particular, the result holds for $\nu=\pi_\thv$).
\end{theorem}
\begin{pf}
Note that $D_n\eqdef\PEs_\thv[\log\pstat(Y_0^{n+1};\thv
)]-
\PEs_\thv[\log\pstat(Y_0^{n};\thv)]$ is a nondecreasing
sequence by \cite{barron1985}, page 1292, and Lemma
\ref{lem:finite-entropy-for-all-n}. Therefore
(\ref{eq:mean-convergence-likelihood-mu}) follows immediately.
As $Y_0^\infty$ is stationary and ergodic under $\PPs_{\thv}$ by
(A1), we can estimate
\begin{eqnarray*}
-\infty&<&D_0\le\sup_{n\ge0}D_n = \ell(\thv) \eqdef
\lim_{n \to\infty}
\PEs_\thv[ n^{-1} \log\pstat(Y_0^n;\thv) ]
\\
&\le&\PEs_\thv\Bigl[
\sup_{x\in\Xset} (\log g_\thv(x,Y_0))^+
\Bigr]<\infty,
\end{eqnarray*}
where we have used (A2). To proceed, we note that
the generalized Shannon--Breiman--McMillan theorem (\cite{barron1985},
Theorem 1), implies that (\ref{eq:as-convergence-likelihood-mu})
holds for $\nu=\pi_\thv$. Therefore, to prove
(\ref{eq:as-convergence-likelihood-mu}) for arbitrary $\nu$, it suffices
to prove the existence of a random variable $C^\nu$ satisfying
$\PPs_{\thv}(0 < C^\nu< \infty)= 1$ and
%
%
\begin{equation}
\label{eq:convergence-ratio}
\lim_{n \to\infty}
\frac{p^\nu(Y_0^n;\thv)}{
\pstat(Y_0^n;\thv)} = C^\nu
, \qquad\PPs_{\thv}\mbox{-a.s.}
\end{equation}
Let $P_n^\nu\eqdef\PP_\thv^\nu(Y_0^n\in\cdot)$ and $\bar
P_n\eqdef\PPs_{\thv}(Y_0^n\in\cdot)$. Then
$p^\nu(Y_0^n;\thv)/\pstat(Y_0^n;\thv)
=\rmd P_n^\nu/\rmd\bar P_n$, and we find that (\ref{eq:convergence-ratio})
holds with $0<C^\nu=\rmd\PP_\thv^{\nu,Y}/\rmd\PPs_\thv^Y<\infty
$ provided
$\PP_\thv^{\nu,Y}\sim\PPs_\thv^Y$. But the latter was already established
in Lemma \ref{lem:ergodic}.
\end{pf}
\begin{rem}
\label{rem:periodicity}
In the case that $g_\thv(x,y)>0$ but $Q_\thv$ is periodic,
the assumption in the above theorem that the initial probability
measure $\nu$
has mass in each periodic class of $Q_\thv$ cannot be eliminated, as
the following example shows. Let $\Xset=\Yset=\{1,2\}$, and let
$Q_\theta$ be the Markov chain with transition probability matrix
$\mathcal{Q}$ and invariant measure $\pi$ (independent of $\theta$)
\[
\mathcal{Q} = \pmatrix{
0 & 1 \cr1 & 0}
,\qquad
\pi= \pmatrix{
1/2 \cr1/2}
.
\]
Then $Q_\theta$ is positive (Harris) recurrent with period 2.
For each $\theta\in\Theta=[0.5,0.9]$, define the observation density
$g_\theta(x,y)$ (with respect to the counting measure)
\[
g_\theta(x,y) = \theta\mathbh{1}_{y=x} + (1-\theta) \mathbh
{1}_{y\ne x},
\]
and let $\thv=0.7$, for example. Then certainly assumptions
(A1) and (A2) are satisfied.

Now consider $\nu=\delta_1$. Then $\nu$ only has mass in one of the two
periodic classes of $Q_\thv$. We can compute the observation likelihood
as follows:
\begin{eqnarray*}
\log p^\nu(Y_0^{2n};\theta)& =&
\sum_{k=0}^{n}\{\mathbh{1}_{Y_{2k}=1}\log\theta+
\mathbh{1}_{Y_{2k}=2}\log(1-\theta)\} \\
&&{} +
\sum_{k=1}^{n}\{\mathbh{1}_{Y_{2k-1}=2}\log\theta+
\mathbh{1}_{Y_{2k-1}=1}\log(1-\theta)\}.
\end{eqnarray*}
A straightforward computation shows that
\begin{eqnarray*}
&&\lim_{n\to\infty}
(2n)^{-1}\log p^\nu(Y_0^{2n};\theta) \\
&&\qquad=
\{\thv\log\theta+ (1-\thv)\log(1-\theta)\}
\mathbh{1}_{X_0=1} \\
&&\qquad\quad{} +
\{(1-\thv)\log\theta+ \thv\log(1-\theta)\}
\mathbh{1}_{X_0=2},\qquad
\PPs_\thv\mbox{-a.s.}
\end{eqnarray*}
Therefore, $\lim_{n\to\infty}n^{-1}\log p^\nu(Y_0^n;\thv)$ is not even
nonrandom $\PPs_\thv\mbox{-a.s.}$, let alone equal to $\ell(\thv
)$. Thus we see
that Theorem \ref{theo:ShannonBreimanMcMillan} does not hold for such
$\nu$.
Moreover, we can compute directly in this example that
\[
\lim_{n\to\infty}{\hat{\theta}}_{\nu,n} = \thv\mathbh
{1}_{X_0=1} +
0.5\mathbh{1}_{X_0=2},\qquad \PPs_\thv\mbox{-a.s.},
\]
so that evidently the maximum likelihood estimator is not consistent when
we choose the initial measure $\nu$. This shows that also
in Theorem \ref{thm:consistency} the assumption that $\nu$
has mass in each periodic class of $Q_\thv$ cannot be eliminated.
\end{rem}

\subsection{Identifiability}
\label{sec:identifiability}

In this section, we establish the identifiability of the parameter
$\theta$. The key issue in the proof consists in showing that the
relative entropy rate between $\pstat(\cdot;\thv)$ and
$p^\lambda(\cdot;\theta)$ may be zero only if $\thv\sim\theta$. Our
proof is based on a very simple and intuitive information-theoretic
device, given as Lemma \ref{lem:general-KL} below, which avoids the
need for an explicit representation of the asymptotic contrast function
as in previous proofs of identifiability.
\begin{defi}
For each $n$, let $\PP_n$ and $\QQ_n$ be probability measures on a
measurable space $(\Zset_n,\Zsigma_n)$. Then $(\QQ_n)$ is
\textit{exponentially separated} from $(\PP_n)$, denoted as
$(\QQ_n)\dashv(\PP_n)$, if there exists a sequence $(A_n)$ of sets
$A_n\in\Zsigma_n$ such that
\begin{eqnarray*}
\liminf_{n\to\infty}\PP_n(A_n)&>&0,\\
\limsup_{n\to\infty}n^{-1}\log\QQ_n(A_n)&<&0 .
\end{eqnarray*}
If $\PP$ and $\QQ$ are probability measures on $(\Yset^{\nset},
\Ysigma^{\otimes\nset})$, then we will write $\QQ\dashv\PP$
if $(\QQ_n)\dashv(\PP_n)$ with $\QQ_n=\QQ(Y_0^n\in\cdot)$ and
$\PP_n=\PP(Y_0^n\in\cdot)$.
\end{defi}
\begin{lem}
\label{lem:general-KL}
If $(\QQ_n)\dashv(\PP_n)$, then $\liminf_{n \to\infty} n^{-1}
\mathrm{KL}(\PP_n||\QQ_n)>0$.
\end{lem}
\begin{pf}
A standard property of the relative entropy
(\cite{dupuisellis1997}, Lemma 1.4.3(g)), states that
for any pair of probability measures $\PP,\QQ$ and measurable set $A$
\[
\mathrm{KL}(\PP||\QQ) \ge\PP(A)\log\PP(A) -
\PP(A)\log\QQ(A)-1,
\]
where $0\log0=0$ by convention. As $x\log x\ge-\rme^{-1}$, we have
\[
\liminf_{n\to\infty}n^{-1}\mathrm{KL}(\PP_n||\QQ_n)
\ge\Bigl(\liminf_{n\to\infty}\PP(A_n)\Bigr)\Bigl(
-\limsup_{n\to\infty}n^{-1}\log\QQ(A_n)\Bigr).
\]
The result follows directly.
\end{pf}

As a consequence of this result, we obtain positive entropy rates:
%
%
\begin{equation}
\label{eq:entropy-rate-simple}
\PP_\theta^{\nu,Y}\dashv\PPs_\thv^Y
\quad\Longrightarrow\quad
\liminf_{n \to\infty}
\PEs_{\thv} \biggl[ n^{-1}
\log\frac{\pstat(Y_0^n;\thv)}{
p^\nu(Y_0^n;\theta)} \biggr]
> 0 .
\end{equation}
This yields identifiability of the asymptotic contrast function in a
very simple and natural manner. It turns out that the exponential
separation assumption $\PP_\theta^{\nu,Y}\dashv\PPs_\thv^Y$
always holds
when the Markov chain $\PP_\theta^\nu$ is $V$-uniformly ergodic and
$\nu(V)<\infty$; this is proved in Section \ref{sec:v-uniform-separation}
below. This observation allows us to establish the consistency of the MLE
in a large class of models.

There is an additional complication that arises in our proof of
consistency. Rather than (\ref{eq:entropy-rate-simple}),
the following result turns out to be of crucial importance:
\[
\liminf_{n \to\infty} \PEs_{\thv} \biggl[
n^{-1} \log\frac{\pstat(Y_0^n;\thv)}{
p^\lambda(Y_0^n;\theta)} \biggr] > 0 .
\]
This result seems almost identical to (\ref{eq:entropy-rate-simple}).
However, note that the probability measure $\nu$ is replaced here by
$\lambda$, the dominating measure on $(\Xset,\Xsigma)$, which may only
be $\sigma$-finite [$\lambda(\Xset)=\infty$]. In this case, a direct
application of Lemma \ref{lem:general-KL} is not possible since
${y}_{0}^{n} \mapsto p^\lambda({y}_{0}^{n};\theta)$ is not a
probability density:
\[
\int p^\lambda({y}_{0}^{n};\theta) \,\rmd{y}_{0}^{n}
=\lambda(\Xset)=\infty.
\]
Nevertheless, the following lemma allows us to reduce the proof in the
case of an improper initial measure $\lambda$ to an application of Lemma
\ref{lem:general-KL}.
\begin{lem}
\label{lem:ptilde}
Assume \textup{(A4)}. For $\theta\not\sim\thv$
such that $p^\lambda({Y}_{0}^{r_\theta};\theta)>0$ $\PPs_\thv\mbox
{-a.s.}$,
there exists a probability measure
$\PPt^\lambda_\theta$ on $(\Yset^\nset,\Ysigma^{\otimes\nset})$
such that
\[
\PPt^\lambda_\theta(Y_0^n\in A) = \int\mathbh{1}_A(y_0^n)
\frac{ p^\lambda({y}_{0}^{n};\theta)}{
p^\lambda({y}_{0}^{r_\theta};\theta)}
\pstat({y}_{0}^{r_\theta};\thv)\,
\rmd y_0^n
\]
for all $n\ge r_\theta$ and $A\in\Ysigma^{\otimes(n+1)}$.
\end{lem}
\begin{pf}
As $\int p^\lambda({y}_{0}^{n};\theta) \,\rmd y_{r_\theta+1}^n =
p^\lambda({y}_{0}^{r_\theta};\theta)<\infty$ $\PPs_\thv$-a.e.
for all $n\geq r_\theta$ by Fubini's theorem and assumption
(A4), we can define for $n\geq r_\theta$
%
%
\begin{equation}
\label{eq:definition-h}
\pt^\lambda({y}_{0}^{n},\theta) \eqdef
p^\lambda({y}_{0}^{n};\theta)
\frac{\pstat({y}_{0}^{r_\theta};\thv)}{
p^\lambda({y}_{0}^{r_\theta};\theta)}
<\infty,\qquad\rmd y_0^n\mbox{-a.e.}
\end{equation}
Note that, by construction, $\{\pt^\lambda(y_0^n,\theta) \,\rmd
y_0^n\dvtx
n\ge r_\theta\}$ is a consistent family of probability measures. By the
extension theorem, we may construct a probability measure
$\PPt^\lambda_\theta$ on $(\Yset^\nset,\Ysigma^{\otimes\nset})$
such that $\PPt^\lambda_\theta(Y_0^n\in A) = \int\mathbh
{1}_A(y_0^n)
\pt^\lambda(y_0^n,\theta) \,\rmd y_0^n$.
\end{pf}
\begin{theorem}
\label{thm:identifiability}
Assume \textup{(A2)}, \textup{(A4)} and
\textup{(A6)}. Then for every $\theta\not\sim\thv$
%
%
\begin{equation}
\label{eq:key-identifiability-0}
\liminf_{n \to\infty} \PEs_{\thv} \biggl[
n^{-1} \log\frac{\pstat(Y_0^n;\thv)}{
p^\lambda(Y_0^n;\theta)} \biggr] > 0.
\end{equation}
\end{theorem}
\begin{pf}
Fix $\theta\not\sim\thv$. Let us assume first that
$\PPs_\thv(p^\lambda({Y}_{0}^{r_\theta};\theta)=0)>0$. As we have
$\int p^\lambda({y}_{0}^{n};\theta) \,\rmd y_{r_\theta+1}^n =
p^\lambda({y}_{0}^{r_\theta};\theta)$ by Fubini's theorem, it must be
the case that $\PPs_\thv(p^\lambda(Y_0^n;\theta)=0)>0$ for all
$n\geq r_\theta$, so that the expression in
(\ref{eq:key-identifiability-0}) is clearly equal to $+\infty$.
Therefore, in this case, the claim holds trivially.

We may therefore assume that $p^\lambda({Y}_{0}^{r_\theta};\theta)>0$
$\PPs_\thv\mbox{-a.s.}$ Let $\pt^\lambda({y}_{0}^{n},\theta)$ be as
in the proof of Lemma \ref{lem:ptilde}. Note that
$\PEs_\thv[|{\log\pstat(Y_0^{r_\theta};\thv)}|]<\infty$ by Lemma
\ref{lem:finite-entropy-for-all-n}, while $\PEs_\thv[(\log
p^\lambda({Y}_{0}^{r_\theta};\theta))^+]<\infty$ by assumption
(A4). Then we find
\[
\liminf_{n \to\infty} \PEs_{\thv} \biggl[
n^{-1} \log\frac{\pstat(Y_0^n;\thv)}{
\pt^\lambda(Y_0^n;\theta)} \biggr] \le
\liminf_{n \to\infty} \PEs_{\thv} \biggl[
n^{-1} \log\frac{\pstat(Y_0^n;\thv)}{
p^\lambda(Y_0^n;\theta)} \biggr].
\]
Assumption (A6) gives $\PPt^\lambda_\theta\dashv
\PPs_\thv^Y$. Therefore, (\ref{eq:key-identifiability-0}) follows
from Lemma \ref{lem:general-KL}.
\end{pf}

\subsection{Consistency of the MLE}
\label{sec:consistency}

Proofs of convergence of the MLE\break \mbox{typically} require to establish the
convergence of the normalized likelihood $n^{-1}\log
p^\nu(Y_0^n;\theta)$ $\PPs_\thv\mbox{-a.s.}$ for any parameter
$\theta$. The
existence of a limit follows from the Shannon--Breiman--McMillan theorem
when $\theta=\thv$ (as in Theorem \ref{theo:ShannonBreimanMcMillan}),
but is far from clear for other $\theta$. In \cite{leroux1992}, the
convergence of $n^{-1}\log p^\nu(Y_0^n;\theta)$ is established using
Kingman's subadditive ergodic theorem. This approach fails in the
present setting, as $\log p^\nu(Y_0^n;\theta)$ may not be subadditive
even up to a constant.

The approach adopted here is inspired by \cite{leroux1992}. We note,
however, that it is not necessary to prove convergence of $n^{-1}\log
p^\nu(Y_0^n;\theta)$ as long as it is asymptotically bounded away from
$\ell(\thv)$, the likelihood of the true parameter. It therefore suffices
to bound $n^{-1}\log p^\nu(Y_0^n;\theta)$ above by an auxiliary sequence
that is bounded away from $\ell(\thv)$. Here the asymptotics of
$n^{-1}\log p^\lambda(Y_0^n;\theta)$ come into play.
\begin{lem}
\label{lem:locally-bounded}
Assume \textup{(A1)--(A6)}.
Then, for any $\theta\not\sim\thv$, there exists an integer
$n_\theta$ and $\eta_\theta>0$ such that
$B(\theta, \eta_\theta)\subseteq\mathcal{U}_\theta$ and
\begin{eqnarray*}
&&\frac{1}{n_\theta+l}
\PEs_\thv\Bigl[
\sup_{\thetapr\in B(\theta,\eta_\theta)}
\log p^\lambda(Y_{0}^{n_\theta};\thetapr)
\Bigr]\\
&&\qquad{}
+ {\frac{1}{n_\theta+l}
\sup_{\thetapr\in B(\theta,\eta_\theta)}
\log}|q_{\thetapr}|_\infty\\
&&\qquad{}+ \frac{l-1}{n_\theta+l}
\PEs_\thv\Bigl[
\sup_{\thetapr\in B(\theta,\eta_\theta)}
\sup_{x\in\Xset} (\log g_\thetapr(x,Y_0))^+
\Bigr]
<\ell(\thv).
\end{eqnarray*}
Here $B(\theta,\eta)\subseteq\Theta$ is the ball of radius $\eta>0$
centered at $\theta\in\Theta$.
\end{lem}
\begin{pf}
By (\ref{eq:mean-convergence-likelihood-mu}) and Theorem
\ref{thm:identifiability}, $\limsup_n n^{-1}\PEs_\thv[ \log
p^\lambda(Y_0^n;\theta)]<\ell(\thv)$.
Using (A4), this implies that there exists a (nonrandom)
integer $n_\theta>r_\theta$ such that
%
%
\begin{eqnarray}\label{eq:firstBoundN_theta}
&&\frac{1}{n_\theta+l}
\PEs_\thv[
\log p^\lambda(Y_{0}^{n_\theta};\theta)
]
+ {\frac{1}{n_\theta+l}
\sup_{\thetapr\in\mathcal{U}_\theta}
\log}|q_{\thetapr}|_\infty\nonumber\\[-8pt]\\[-8pt]
&&\qquad{}+ \frac{l-1}{n_\theta+l}
\PEs_\thv\Bigl[
\sup_{\thetapr\in\mathcal{U}_\theta}
\sup_{x\in\Xset} (\log g_\thetapr(x,Y_0))^+
\Bigr]
<\ell(\thv).\nonumber
\end{eqnarray}
For any $\eta> 0$ such that $B(\theta,\eta) \subseteq\calU_\theta$,
we have
\begin{eqnarray*}
&&
\sup_{\thetapr\in B(\theta, \eta)}
\log p^\lambda({Y}_{0}^{n_\theta}; \thetapr)\\
&&\qquad\leq
\sup_{\thetapr\in\calU_\theta} (\log
p^\lambda({Y}_{0}^{r_\theta}; \thetapr))^+
\\
&&\qquad\quad{}+ \sum_{k=r_\theta+1}^{n_\theta}
\sup_{\thetapr\in\calU_\theta}\sup_{x\in\Xset}
(\log g_\thetapr(x,Y_k))^+,
\end{eqnarray*}
where the right-hand side does not depend on $\eta$ and is integrable.
But then
%
%
\begin{eqnarray}\label{eq:techn1}
&&\limsup_{\eta\downarrow0}
\PEs_\thv\Bigl[ \sup_{\thetapr\in B(\theta, \eta)}
\log p^\lambda({Y}_{0}^{n_{\theta}}; \thetapr)
\Bigr]\nonumber\\
&&\qquad\leq
\PEs_\thv\Bigl[
\limsup_{\eta\downarrow0} \sup_{\thetapr\in B(\theta, \eta)}
\log p^\lambda({Y}_{0}^{n_\theta} ; \thetapr)\Bigr]\\
&&\qquad \le
\PEs_\thv[
\log p^\lambda({Y}_{0}^{n_\theta} ; \theta)]
,\nonumber
\end{eqnarray}
by (A5) and Fatou's lemma.
Together (\ref{eq:techn1}) and (\ref{eq:firstBoundN_theta}) complete
the proof.
\end{pf}
\begin{pf*}{Proof of Theorem \ref{thm:consistency}}
Since, by Theorem \ref{theo:ShannonBreimanMcMillan},
$\lim_{n \rightarrow\infty} n^{-1}\ell_{\nu,n}(\thv) = \ell(\thv)$,
$\PPs_\thv$-a.s., it is sufficient to prove that for any closed set
$\Cset\subset\Theta$ such that \mbox{$\Cset\cap[\thv]=\varnothing$}
\[
\limsup_{n \to\infty} \sup_{\thetapr\in\Cset}
n^{-1}\ell_{\nu,n}(\thetapr) < \ell(\thv) ,\qquad
\PPs_\thv\mbox{-a.s.}
\]
Now note that $\{B(\theta,\eta_\theta)\dvtx\theta\in\Cset\}$ is a
cover of
$\Cset$, where $\eta_\theta$ are defined in Lem\-ma~\ref{lem:locally-bounded}.
As $\Theta$ is compact, $\Cset$ is also compact and thus admits a finite
subcover $\{B(\theta_i,\eta_{\theta_i})\dvtx\theta_i\in\Cset
,i=1,\ldots,N\}$.
It therefore suffices to show that
\[
\limsup_{n \to\infty} \sup_{\thetapr\in
B(\theta,\eta_\theta)\cap\Cset}
n^{-1}\ell_{\nu,n}(\thetapr) < \ell(\thv) ,\qquad
\PPs_\thv\mbox{-a.s.}
\]
for any $\theta\not\sim\thv$. Fix $\theta\not\sim\thv$ and let
$\eta_\theta$ and $n_\theta$ be as in Lemma \ref{lem:locally-bounded}.
Note that
%
%
\begin{eqnarray}
\label{eq:splitlik:nu}
p^\nu(y_0^n;\thetapr) &\le&
p^\nu(y_0^m;\thetapr) p^\lambda(y_{m+l}^n;\thetapr)
g_\thetapr^*(y_{m+1}^{m+l-1}) |q_{\thetapr}|_\infty,
\\
\label{eq:splitlik:lambda}
p^\lambda(y_j^n;\thetapr) &\le&
p^\lambda(y_j^m;\thetapr) p^\lambda(y_{m+l}^n;\thetapr)
g_\thetapr^*(y_{m+1}^{m+l-1}) |q_{\thetapr}|_\infty,
\end{eqnarray}
for any $j\le m$, $m+l\le n$ and $\thetapr\not\sim\thv$, where
$g_\theta^*(y_i^j) \eqdef\prod_{\ell=i}^j\sup_{x\in\Xset
}g_\theta(x,y_\ell)$.
We can therefore estimate, for all $n$ sufficiently large,
\begin{eqnarray*}
\ell_{\nu,n}(\thetapr)
&\le&
\frac{1}{n_\theta+l}
\sum_{r=1}^{n_\theta+l}\{
\ell_{\nu,r-1}(\thetapr)
+ \log p^\lambda(Y_{l+r-1}^n;\thetapr)
+ \log(g_\thetapr^*(Y_r^{l+r-2}) |q_{\thetapr}|_\infty)
\}
\\
&\le&
\frac{1}{n_\theta+l}
\sum_{r=1}^{n_\theta+l}
\sum_{k=1}^{i(n)-1}
\bigl\{
\log p^\lambda\bigl(Y_{(n_\theta+l)(k-1)+l+r-1}^{(n_\theta+l)k
+r-1};\thetapr\bigr)
+ {\log}|q_{\thetapr}|_\infty\bigr\}
\\
&&{}
+ \frac{1}{n_\theta+l} \sum_{r=1}^{n_\theta+l}\sum_{k=1}^{i(n)-1}
\log g_\thetapr^*\bigl(Y_{(n_\theta+l)(k-1)+r}^{
(n_\theta+l)(k-1)+l+r-2}\bigr)
\\
&&{}
+ \frac{1}{n_\theta+l} \sum_{r=1}^{n_\theta+l}
\log\bigl(g_\thetapr^*(Y_0^{r-1})
g_\thetapr^*\bigl(Y_{(n_\theta+l)(i(n)-1)+r}^n\bigr)\bigr)
\\ &=&
\frac{1}{n_\theta+l}
\sum_{r=1}^{(n_\theta+l)(i(n)-1)}
\Biggl\{
\log p^\lambda(Y_{l+r-1}^{n_\theta+l+r-1};\thetapr)
+ \sum_{k=0}^{l-2} \sup_{x\in\Xset} \log g_\thetapr(x,Y_{k+r})
\Biggr\}
\\
&&{}
+ {\bigl(i(n)-1\bigr)\log}|q_{\thetapr}|_\infty
+ \frac{1}{n_\theta+l} \sum_{r=1}^{n_\theta+l}
\log g_\thetapr^*(Y_0^{r-1})
\\
&&{}
+ \frac{1}{n_\theta+l} \sum_{r=1}^{n_\theta+l}
\sum_{k=(n_\theta+l)(i(n)-1)+r}^n
\sup_{x\in\Xset}
\log g_\thetapr(x,Y_k)
,
\end{eqnarray*}
where $i(n)\eqdef\max\{m\in\mathbb{N}\dvtx m(n_\theta+l)\le n\}$.
Here we have applied (\ref{eq:splitlik:nu}) with $m=r-1$ in the first
inequality, while we have repeatedly applied (\ref{eq:splitlik:lambda})
for every $m=(n_\theta+l)k+r-1$, $k\le i(n)-2$ in the second inequality,
together with the simple estimates
$\ell_{\nu,r-1}(\theta')\le\log g_{\theta'}^*(Y_0^{r-1})$
and
\begin{eqnarray*}
&&p^\lambda\bigl(Y_{(n_\theta+l)(i(n)-2)+l+r-1}^n;\theta'\bigr)\\
&&\qquad\le
p^\lambda\bigl(Y_{(n_\theta+l)(i(n)-2)+l+r-1}^{(n_\theta+l)(i(n)-1)+r-1};
\theta'\bigr) g^*_{\theta'}\bigl(Y_{(n_\theta+l)(i(n)-1)+r}^n
\bigr).
\end{eqnarray*}
We can now estimate, for all $n$ sufficiently large,
\begin{eqnarray*}
&&\sup_{\thetapr\in B(\theta,\eta_\theta)\cap\Cset}
\ell_{\nu,n}(\thetapr) \\
&&\qquad\le
\frac{1}{n_\theta+l}
\sum_{r=1}^{(n_\theta+l)(i(n)-1)}
\sup_{\thetapr\in B(\theta,\eta_\theta)\cap\Cset}
\log p^\lambda(Y_{l+r-1}^{n_\theta+l+r-1};\thetapr)
\\
&&\qquad\quad{}
+ \sum_{k=0}^{l-2}
\frac{1}{n_\theta+l}
\sum_{r=1}^{(n_\theta+l)(i(n)-1)}
\sup_{\thetapr\in B(\theta,\eta_\theta)\cap\Cset}
\sup_{x\in\Xset} (\log g_\thetapr(x,Y_{k+r}))^+
\\
&&\qquad\quad{}
+ {\bigl(i(n)-1\bigr)\sup_{\thetapr\in B(\theta,\eta_\theta)\cap\Cset}
\log}|q_{\thetapr}|_\infty\\
&&\qquad\quad{}
+ \frac{1}{n_\theta+l} \sum_{r=1}^{n_\theta+l}
\sup_{\thetapr\in B(\theta,\eta_\theta)\cap\Cset}
\log g_\thetapr^*(Y_0^{r-1})
\\
&&\qquad\quad{}
+ \sum_{k=n-2(n_\theta+l)+1}^n
\sup_{\thetapr\in B(\theta,\eta_\theta)\cap\Cset}
\sup_{x\in\Xset}
(\log g_\thetapr(x,Y_k))^+,
\end{eqnarray*}
where we have used that $(n_\theta+l)(i(n)-1)+r\ge n-2(n_\theta+l)+1$ to
estimate the last term. But as $i(n)/n\to(n_\theta+l)^{-1}$ as $n\to
\infty$,
we find that
\begin{eqnarray*}
\limsup_{n \to\infty} \sup_{\thetapr\in B(\theta,\eta_\theta)
\cap\Cset}n^{-1} \ell_{\nu,n}(\thetapr)
&\le&
\frac{1}{n_\theta+l}
\PEs_\thv\Bigl[
\sup_{\thetapr\in B(\theta,\eta_\theta)}
\log p^\lambda(Y_{0}^{n_\theta};\thetapr)
\Bigr]
\\
&&{}
+ \frac{l-1}{n_\theta+l}
\PEs_\thv\Bigl[
\sup_{\thetapr\in B(\theta,\eta_\theta)}
\sup_{x\in\Xset} (\log g_\thetapr(x,Y_0))^+
\Bigr]\\
&&{}+ {\frac{1}{n_\theta+l}
\sup_{\thetapr\in B(\theta,\eta_\theta)}
\log}|q_{\thetapr}|_\infty
\\
&<& \ell(\thv)
\end{eqnarray*}
by (A4), Birkhoff's ergodic theorem,
Lemma \ref{lem:locally-bounded}, and the elementary fact that
$\lim_{n}\frac{1}{n}\sum_{k=n-r+1}^n\xi_k=\lim_{n}
\frac{1}{n}\sum_{k=1}^n\xi_k-\lim_{n}
\frac{1}{n}\sum_{k=1}^{n-r}\xi_k=0$
for any stationary ergodic sequence $(\xi_k)_{k\ge0}$ with
$\mathbb{E}(|\xi_1|)<\infty$. This completes the proof.
\end{pf*}

\section{Exponential separation and $V$-uniform ergodicity}
\label{sec:v-uniform-separation}

As is explained in Remark \ref{rem:sketch-proof-v}, the key step in
establishing assumption (A6) is to obtain a type of large
deviations property. The following Azuma--Hoeffding type inequality
provides what is needed in the $V$-uniformly ergodic case.
\begin{theorem}
\label{thm:azuma-hoeffding-v-unif}
Assume that $Q_\theta$ is $V_\theta$-uniformly ergodic. Fix $s \geq0$,
and let $f\dvtx \Yset^{s+1} \to\rset$ be such that $|f|_\infty<\infty$.
Then there exists a constant $K$ such that
\[
\PP_\theta^\nu\Biggl( \Biggl|\sum_{i=1}^n
\{f({Y}_{i}^{i+s}) - \PEs_\theta[f({Y}_{0}^{s})] \}
\Biggr| \geq t \Biggr) \leq
K \nu(V) \exp\biggl[ - \frac{1}{K}
\biggl( \frac{t^2}{n} \wedge t \biggr) \biggr]
\]
for any probability measure $\nu$ and any $t > 0$.
\end{theorem}

We will first use this result in Section \ref{sec:proof-v-uniform} to
prove Theorem \ref{thm:v-uniformly-ergodic}. In Section
\ref{sec:exponential-inequality}, we will establish a general
Azuma--Hoeffding type large deviations inequality for $V$-uniformly
ergodic Markov chains, which forms the basis for the proof of
Theorem~\ref{thm:azuma-hoeffding-v-unif}. Finally, Section
\ref{sec:proof-azuma-hoeffding-v-unif} completes the proof of Theorem~\ref{thm:azuma-hoeffding-v-unif}.

\subsection{\texorpdfstring{Proof of Theorem
\protect\ref{thm:v-uniformly-ergodic}}{Proof of Theorem 2}}
\label{sec:proof-v-uniform}

We begin by proving that exponential separation holds under the
$V$-uniform ergodicity assumption.
\begin{prop} \label{prop:kullBpropre}
Assume \textup{(A1)} and \textup{(A6$'$)}. For any
$\theta\not\sim\thv$ with
$p^\lambda({Y}_{0}^{r_\theta};\theta)>0$ $\PPs_\thv\mbox{-a.s.}$
and probability measure $\nu$ such that
$\nu(V_\theta)<\infty$, we have $\PP_\theta^{\nu,Y}\dashv\PPs
_\thv^Y$.
\end{prop}
\begin{pf}
Fix $\theta\not\sim\thv$. As $\PPs_\theta^Y\ne\PPs_\thv^Y$ by
assumption
(A6$'$), there exists an integer $s\geq0$ and a bounded
measurable function $h\dvtx \Yset^{s+1} \to\rset$ such that
$\PEs_\theta[h({Y}_{0}^{s})]=0$ and
$\PEs_\thv[h({Y}_{0}^{s}) ]=1$. Define for $n\ge s$ the
set $A_n\in\Ysigma^{\otimes(n+1)}$ as
\[
A_n \eqdef
\Biggl\{{y}_{0}^{n} \in\Yset^{n+1}\dvtx
\Biggl| \frac{1}{n-s}\sum_{i=1}^{n-s}
h({y}_{i}^{i+s}) \Biggr| \geq\frac{1}{2} \Biggr\}.
\]
As $Y_0^\infty$ is stationary and ergodic under $\PPs_\thv$ by
(A1), Birkhoff's ergodic theorem gives
$\PPs_\thv^Y(A_n)\to1$ as $n\to\infty$. On the other hand,
Theorem \ref{thm:azuma-hoeffding-v-unif} shows that
$\limsup_{n\to\infty}n^{-1}\log\PP_\theta^{\nu,Y}(A_n)<0$.
Thus, we have established $\PP_\theta^{\nu,Y}\dashv\PPs_\thv^Y$.
\end{pf}

Proposition \ref{prop:kullBpropre} is not sufficient to establish
(A6), however: the problem is that we are interested in the
case where $\nu$ is not a probability measure, but the $\sigma$-finite
measure $\lambda$. What remains is to reduce this problem to an
application of Proposition~\ref{prop:kullBpropre}. To this end, we will
use the following lemma.
\begin{lem}
\label{lem:lambda-b-theta}
Assume \textup{(A4)}, and fix
$\theta\not\sim\thv$ such that
$p^\lambda({Y}_{0}^{r_\theta};\theta)>0$ $\PPs_\thv\mbox{-a.s.}$
For any $B\in\Ysigma^{\otimes(r_\theta+1)}$
such that $\PPs_\thv(Y_0^{r_\theta}\in B)>0$, define the probability
measure
\[
\lambda_{B,\theta}(A) =
\PEs_\thv\biggl(
\int\mathbh{1}_A(x_{r_\theta+1}) \frac{p^\lambda(\rmd x_{r_\theta+1},
Y_0^{r_\theta};\theta)}{p^\lambda(Y_0^{r_\theta};\theta)}
\Big|
Y_0^{r_\theta}\in B
\biggr)
\]
on $(\Xset,\Xsigma)$. Then we have
\[
\PPt^\lambda_\theta(Y_0^{r_\theta}\in B,
Y_{r_\theta+1}^n\in A) = \PPs_\thv(Y_0^{r_\theta}\in B)
\PP_\theta^{\lambda_{B,\theta}}(Y_0^{n-r_\theta-1}\in A)
\]
for any set $A\in\Ysigma^{\otimes(n-r_\theta)}$.
\end{lem}
\begin{pf}
Note that by assumption (A4), $\PPt^\lambda_\theta$ is well
defined (as shown in Lem\-ma~\ref{lem:ptilde}) and
$0<p^\lambda(Y_0^{r_\theta};\theta)<\infty$ $\PPs_\thv\mbox{-a.s.}$
Moreover, as $p^\lambda(y_0^{r_\theta};\theta) =
\int p^\lambda(\rmd x_{r_\theta+1}$, $y_0^{r_\theta};\theta)$,
we find that $\lambda_{B,\theta}$ is indeed a probability measure
on $(\Xset,\Xsigma)$.

Let $B\in\Ysigma^{\otimes(r_\theta+1)}$ be such that
$\PPs_\thv(Y_0^{r_\theta}\in B)>0$. Then for any $n>r_\theta$
\begin{eqnarray*}
&&\PPt^\lambda_\theta(Y_0^{r_\theta}\in B,
Y_{r_\theta+1}^n\in A)\\
&&\qquad=
\int\mathbh{1}_A(y_{r_\theta+1}^n) \mathbh{1}_B(y_0^{r_\theta})
p^\lambda({y}_{0}^{n};\theta)
\frac{\pstat({y}_{0}^{r_\theta};\thv)}{
p^\lambda({y}_{0}^{r_\theta};\theta)}\,
\rmd y_0^n \\
&&\qquad=
\PPs_\thv(Y_0^{r_\theta}\in B)
\int\biggl[\int
\mathbh{1}_A(y_{r_\theta+1}^n)
p^{x_{r_\theta+1}}(y_{r_\theta+1}^n;\theta) \,\rmd y_{r_\theta+1}^n\biggr]
\lambda_{B,\theta}(\rmd x_{r_\theta+1})\\
&&\qquad=
\PPs_\thv(Y_0^{r_\theta}\in B)
\PP_\theta^{\lambda_{B,\theta}}(Y_0^{n-r_\theta-1}\in A),
\end{eqnarray*}
where we used $p^\lambda({y}_{0}^{n};\theta)=\int
p^\lambda(\rmd x_{m+1},y_0^m;\theta)
p^{x_{m+1}}(y_{m+1}^n;\theta)$ for $n>m$.
\end{pf}

We can now complete the proof of Theorem \ref{thm:v-uniformly-ergodic}.
\begin{pf*}{Proof of Theorem \ref{thm:v-uniformly-ergodic}}
Fix $\theta\not\sim\thv$ such that
$p^\lambda({Y}_{0}^{r_\theta};\theta)>0$ $\PPs_\thv\mbox{-a.s.}$, and
define
\[
B = \biggl\{
y_0^{r_\theta}\dvtx
\int V_\theta(x_{r_\theta+1})
\frac{p^\lambda(\rmd x_{r_\theta+1},{y}_{0}^{r_\theta};\theta)}
{p^\lambda({y}_{0}^{r_\theta};\theta)}
\le K\biggr\}.
\]
By (A6$'$), we can choose $K$ sufficiently large so that
$\PPs_\thv(Y_0^{r_\theta}\in B)>0$. Consequently
$\lambda_{B,\theta}(V_\theta)\le K <\infty$ by construction.
As in the proof of Proposition \ref{prop:kullBpropre}, it follows that
there exists a sequence of sets $A_n\in\Ysigma^{\otimes(n-r_\theta)}$
such that
\[
\lim_{n\to\infty}\PPs_\thv(Y_0^{n-r_\theta-1}\in A_n)=1,\qquad
\limsup_{n\to\infty}n^{-1}
\log\PP_\theta^{\lambda_{B,\theta}}(Y_0^{n-r_\theta-1}\in A_n)
<0.
\]
Define the sets
\[
\tilde A_n \eqdef\{y_0^n\dvtx y_0^{r_\theta}\in B,
y_{r_\theta+1}^n\in A_n\}.
\]
Using the stationarity of $\PPs_\thv$ and
Lemma \ref{lem:lambda-b-theta}, it follows that
\[
\lim_{n\to\infty}\PPs_\thv(Y_0^n\in\tilde A_n) =
\PPs_\thv(Y_0^{r_\theta}\in B)>0,\qquad
\limsup_{n\to\infty}
n^{-1}\log\PPt_\theta^\lambda(Y_0^n\in\tilde A_n)<0.
\]
This establishes (A6).
\end{pf*}

Finally, let us prove Proposition \ref{prop:mixing-for-all-initial}.
\begin{pf*}{Proof of Proposition \ref{prop:mixing-for-all-initial}}
Fix $\theta\not\sim\thv$ such that
$p^\lambda({Y}_{0}^{r_\theta};\theta)>0$ $\PPs_\thv\mbox{-a.s.}$, and
let $B=\Yset^{r_\theta+1}$. By (A6$''$), there exists a
sequence of sets $A_n\in\Ysigma^{\otimes(n-r_\theta)}$ such that
\[
\liminf_{n\to\infty}\PPs_\thv(Y_0^{n-r_\theta-1}\in A_n)>0,\qquad
\limsup_{n\to\infty}n^{-1}
\log\PP_\theta^{\lambda_{B,\theta}}(Y_0^{n-r_\theta-1}\in A_n)
<0.
\]
Assumption (A6) now follows easily from the stationarity of
$\PPs_\thv$ and Lemma \ref{lem:lambda-b-theta}.
\end{pf*}

\subsection{An Azuma--Hoeffding inequality}
\label{sec:exponential-inequality}

This section is somewhat independent of the remainder of the paper. We
will prove a general Azuma--Hoeffding type large deviations inequality
for $V$-uniformly ergodic Markov chains, on which the proof of Theorem
\ref{thm:azuma-hoeffding-v-unif} will be based (see Section
\ref{sec:proof-azuma-hoeffding-v-unif}). The following result may
be seen as an extension of the Azuma--Hoeffding inequality obtained in
\cite{glynnormoneit2002} for uniformly ergodic Markov chains, and the
proof of our result is similar to the proof of the Bernstein-type
inequality in \cite{adamczak2008}, Theorem~6.
\begin{theorem}
\label{thm:exponential-inequality}
Let $(X_k)_{k\ge0}$ be a Markov chain in $(\Xset,\Xsigma)$ with transition
kernel $Q$ and initial measure $\eta$ under the probability measure
$\PP^\eta$.
Assume that the transition kernel $Q$ is $V$-uniformly ergodic, and denote
by $\pi$ its unique invariant measure.
Then there exists a constant $K$ such that
\[
\PP^\eta\Biggl( \Biggl|
\sum_{i=1}^n \{f (X_i) - \pi(f) \} \Biggr| \geq t
\Biggr)
\leq
K \eta(V) \exp\biggl[ - \frac{1}{K}
\biggl( \frac{t^2}{n|f|_\infty^2} \wedge\frac{t}{|f|_\infty}
\biggr) \biggr]
\]
for any probability measure $\eta$, bounded function
$f\dvtx\Xset\to\rset$, and $t>0$.
\end{theorem}
\begin{rem}
The exponential bound of Theorem \ref{thm:exponential-inequality} has a
Bernstein-type tail, unlike the usual Azuma--Hoeffding bound. However,
unlike the Bernstein inequality, the tail behavior is determined only by
$|f|_\infty$, and not by the variance of $f$. We therefore still refer
to this inequality as an Azuma--Hoeffding bound. It is shown in
\cite{adamczak2008} by means of a counterexample that $V$-uniformly
ergodic Markov chains do not admit, in general, a Bernstein bound
of the type available for independent random variables (the bound in
\cite{adamczak2008} depends on the variance at the cost of an extra
logarithmic factor, which precludes its use for our purposes).
\end{rem}

Throughout this section, we let $(X_k)_{k\ge0}$ be as in Theorem
\ref{thm:exponential-inequality}. For simplicity, we work with a
generic constant $K$ which may change from line to line.

Before we turn to the proof of Theorem \ref{thm:exponential-inequality},
let us recall some standard facts from the theory of $V$-uniformly
ergodic Markov chains. It is well known (\cite{meyntweedie1993}, Chapter 16),
that $V$-uniform ergodicity in the sense of
Definition \ref{defi:v-uniformly-ergodic} implies (and is essentially
equivalent to) the following properties:

\textit{Minorization condition}. There exist a set $C\in\Xsigma$,
an integer $m$, a probability measure $\nu$ on $(\Xset,\Xsigma)$ and
a constant $\varepsilon> 0$ such that
%
%
\begin{equation}
\label{eq:minorization-condition}
Q^m(x,A) \geq\varepsilon\nu(A)
\qquad\mbox{for all }x \in C\mbox{ and all }A \in\Xsigma.
\end{equation}

\textit{Foster--Lyapunov drift condition}. There exists a measurable
function $V\dvtx
\Xset\mapsto[ 1,\infty)$, $\lambda\in[ 0,1 )$, and $b
< \infty$, such that $\sup_{x \in C} V(x) < \infty$ and
%
%
\begin{equation}
\label{eq:drift-condition}
\Q V(x) \leq\lambda V(x) + b \mathbh{1}_C(x)
\qquad\mbox{for all }x\in\Xset.
\end{equation}

The set $C$ in the minorization condition is referred to as a
\textit{$(\nu,m)$-small set} (see \cite{meyntweedie1993} for
extensive discussion). For future reference, let us note that
\[
1\le\pi(V) = (1-\lambda)^{-1}\pi(QV-\lambda V) \le
(1-\lambda)^{-1}b \pi(C)<\infty,
\]
which shows that $\pi(V)<\infty$ and $\pi(C)>0$. Moreover,
\[
\varepsilon\pi(C) \nu(V)\le\pi(Q^mV) = \pi(V)<\infty,
\]
so that necessarily $\nu(V)<\infty$ also.

The proof of Theorem \ref{thm:exponential-inequality} is based on an
embedding of the Markov chain into a wide sense regenerative process
(\cite{kalashnikov1994}, page 360), known as a \textit{splitting
construction}. Let us recall how this can be done. We will employ
the canonical process $\cX_n \eqdef(\tilde{X}_n,d_n)$ on the enlarged
measure space $(\check\Omega,\mathcal{\check F})$, where
$\check\Omega=(\Xset\times\{0,1\})^{\nset}$ and $\mathcal{\check
F}$ is
the corresponding Borel $\sigma$-field. In words, $\tilde{X}_n$ takes
values in $(\Xset,\Xsigma)$ and $d_n$ is a binary random variable.
Define the following stopping times:
\[
\sigma_0 \eqdef\inf\{n\ge0\dvtx\tilde{X}_n\in C\},\qquad
\sigma_{i+1} \eqdef\inf\{n\ge\sigma_i+m\dvtx\tilde{X}_n\in C\}.
\]
We now construct a probability measure $\cPP^{\eta}$ on
$(\check\Omega,\mathcal{\check F})$ with the following properties
(e.g., by means of the Ionescu--Tulcea theorem):
\begin{eqnarray*}
&& (d_n)_{n\ge0}\mbox{ are i.i.d. }\qquad\mbox{with }\cPP^{\eta
}(d_n=1)=\varepsilon
,
\\
&&\tilde{X}_0\mbox{ is independent from
}(d_n)_{n\ge0}\quad\mbox{and}\quad
\cPP^{\eta}(\tilde{X}_0\in\cdot)=\eta,
\\
&&\cPP^{\eta}(\tilde{X}_{n+1}\in\cdot|\tilde{X}_0^n,d_0^\infty)
= Q(\tilde{X}_n, \cdot)\qquad\mbox{on }
\{n<\sigma_0\}\cup\bigcup_{i\ge0}\{\sigma_i+m\le n<\sigma_{i+1}\}
,
\\
&&\cPP^{\eta}(\tilde{X}_{\sigma_i+1}^{\sigma_i+m}\in\cdot|
\tilde{X}_0^{\sigma_i},d_0^\infty)=\cases{
\displaystyle \int\mathbf{q}^{\tilde{X}_{\sigma_i},x}( \cdot) \nu(\rmd x),
&\quad if $d_{\sigma_i}=1$,\vspace*{2pt}\cr
\displaystyle \int\mathbf{q}^{\tilde{X}_{\sigma_i},x}( \cdot)
R(\tilde{X}_{\sigma_i},\rmd x), &\quad if $d_{\sigma_i}=0$.}
\end{eqnarray*}
Here we defined the transition kernel
$R(x,A)\eqdef(1-\varepsilon)^{-1}\{Q^m(x,A)-\varepsilon\nu(A)\}$ for
$x\in C$, and (using that $\Xset$ is Polish to ensure existence) the
regular conditional probability
$\mathbf{q}^{X_0,X_m}(A)\eqdef\PP^\eta(X_1^m\in A|X_0,X_m)$.

The process $(\cX_n)_{n\ge0}$ is not necessarily Markov. However, it
is easily verified that the law of the process $(\tilde{X}_n)_{n\ge0}$
under $\cPP^\eta$ is the same as the law of $(X_n)_{n\ge0}$ under
$\PP^\eta$, so that our original Markov chain is indeed embedded in this
construction. Moreover, at every time $\sigma_n$ such that additionally
$d_{\sigma_n}=1$, we have by construction that $\tilde{X}_{\sigma_n+m}$
is drawn independently from the distribution $\nu$, that is, the process
regenerates in $m$ steps. Let us define the regeneration times as
\[
\csigma_0 \eqdef\inf\{\sigma_i+m\dvtx i\ge0, d_{\sigma_i}=1\},\qquad
\csigma_{n+1} \eqdef\inf\{\sigma_i+m\dvtx\sigma_i\ge\csigma_{n},
d_{\sigma_i}=1\}.
\]
The regenerations will allow us to split the path of the process into
one-dependent blocks, to which we can apply classical large deviations
bounds for independent random variables. We formalize this as the
following lemma.
\begin{lem}
\label{lem:block-sums}
Define for $i\ge0$ the block sums
\[
\xi_i \stackrel{\mathit{def}}{=}\sum_{k=\csigma_i}^{\csigma_{i+1}-1}
\{f(\tilde{X}_k) - \pi(f)\} .
\]
Then $(\xi_i)_{i\ge0}$ are identically distributed, one-dependent,
and $\cPE^\eta(\xi_0)=0$.
\end{lem}
\begin{pf}
First, we note that $\cPP^{\eta}(\cX_{\csigma_i}^{\csigma_{i+1}-1}
\in\cdot|\cX_0^{\csigma_i-m})=\cPP^\nu(\cX_0^{\csigma
_0-1}\in\cdot)$
for all~$i$. It follows directly that $(\xi_i)_{i\ge0}$ are
identically distributed and one-dependent. Moreover, as $\csigma_i$ is
$\sigma\{X_0^{\csigma_i-m}\}$-measurable, we find that the
inter-regeneration times $(\csigma_{i+1}-\csigma_i)_{i\ge0}$ are
independent. Now note that, by the law of large numbers,
\begin{eqnarray*}
\cPE^\eta(\xi_0) &=&
\lim_{n\to\infty}\frac{1}{n}\sum_{i=0}^{n-1}\xi_i =
\lim_{n\to\infty}\frac{1}{n}\sum_{k=\csigma_0}^{\csigma_n-1}
\{f(\tilde{X}_k)-\pi(f)\}  \\
&=&\lim_{n\to\infty}
\Biggl(\frac{1}{n}\sum_{i=0}^{n-1}\{\csigma_{i+1}-\csigma_i\}\Biggr)
\Biggl(\frac{1}{\csigma_n-\csigma_0}\sum_{k=\csigma_0}^{\csigma_n-1}
\{f(\tilde{X}_k)-\pi(f)\}\Biggr).
\end{eqnarray*}
But $\lim_{n\to\infty}\frac{1}{n}\sum_{i=0}^{n-1}\{\csigma
_{i+1}-\csigma_i\}=
\cPE^\eta(\csigma_1-\csigma_0)<\infty$ by the law of large numbers and
(\ref{eq:bound-regeneration-time}) below, while
$\lim_{n\to\infty}\frac{1}{\csigma_n-\csigma_0}\sum_{k=\csigma
_0}^{\csigma_n-1}
\{f(\tilde{X}_k)-\pi(f)\}=0$ by the ergodic theorem for Markov chains.
This completes the proof.
\end{pf}

In the proof of Theorem \ref{thm:exponential-inequality},
we will need that fact that the inter-regeneration times
$\csigma_0$ and $\csigma_{i+1}-\csigma_i$ possess exponential
moments. We presently establish that this is necessarily the
case, adapting the proof of \cite{robertstweedie1999}, Theorem 2.1.
\begin{prop}
There exists a constant $K$ such that
%
%
\begin{equation}
\label{eq:bound-regeneration-time}
\cPE^{\eta}[ \exp(\csigma_0/K) ]
\leq K \eta(V) \quad\mbox{and}\quad
\cPE^{\eta}[ \exp(\{\csigma_{1} - \csigma_0\}/K)
] \leq K
\end{equation}
for every probability measure $\eta$.
\end{prop}
\begin{pf}
We begin by writing
\[
\{\csigma_0-m=n\} =
\bigcup_{j\ge0}
\{d_{\sigma_0},\ldots,d_{\sigma_{j-1}}=0, d_{\sigma_j}=1, \sigma
_j=n\}.
\]
Using the independence of $d_{\sigma_j}$ from
$d_0,\ldots,d_{\sigma_{j-1}},\sigma_j$, we have
\[
\cPP^\eta(\csigma_0-m=n) =
\sum_{j=0}^\infty
\varepsilon(1-\varepsilon)^j
\cPP^\eta(\sigma_j=n|d_{\sigma_0},\ldots,d_{\sigma_{j-1}}=0).
\]
In particular, we can write
\[
\cPE^{\eta}(\rme^{\csigma_0/K}) =
\rme^{m/K}\sum_{j=0}^\infty\varepsilon(1-\varepsilon)^j
\cPE^\eta(\rme^{\sigma_j/K}|d_{\sigma_0},\ldots,d_{\sigma_{j-1}}=0).
\]
Now note that by construction, we have
\[
\cPE^\eta\bigl(\rme^{\{\sigma_j-\sigma_{j-1}-m\}/K}|\cX_0^{\sigma
_{j-1}}\bigr) =
\cPE^{R(\tilde{X}_{\sigma_{j-1}}, \cdot)}(\rme^{\sigma_0/K})
\qquad\mbox{on }\{d_{\sigma_{j-1}}=0\}.
\]
Define $G(K)\eqdef\sup_{x\in C}\cPE^{R(x, \cdot)}(\rme^{\sigma_0/K})$.
It is now easily established that
\[
\cPE^\eta(\rme^{\sigma_j/K}|d_{\sigma_0},\ldots,d_{\sigma_{j-1}}=0)
\le\rme^{jm/K}G(K)^j \cPE^\eta(\rme^{\sigma_0/K}).
\]
We can therefore estimate
\[
\cPE^{\eta}(\rme^{\csigma_0/K}) \le
\frac{\varepsilon\rme^{m/K} \cPE^\eta(\rme^{\sigma_0/K})}{
1-(1-\varepsilon) \rme^{m/K} G(K)},
\]
provided that $(1-\varepsilon) \rme^{m/K} G(K)<1$.

Now note that it follows from \cite{meyntweedie1993}, Theorem 15.2.5, that
%
%
\begin{equation}
\label{eq:meyn-tweedie-estimate-sigma0}
\cPE^x(\rme^{\sigma_0/K}) \le K\{\lambda V(x)+b\mathbh{1}_C(x)\}
\qquad\mbox{for all }x\in\Xset,
\end{equation}
provided $K$\vspace*{1pt} is chosen sufficiently large.
Therefore, it is easily established that $\cPE^\eta(\rme^{\sigma
_0/K})\le
K \eta(V)$ for $K$ sufficiently large. On the other hand, by Jensen's
inequality, $G(K)\le G(\beta)^{\beta/K}$ for $\beta\le K$. As
$G(\beta)<\infty$ for some $\beta$ by
(\ref{eq:meyn-tweedie-estimate-sigma0}), we have $G(K)\to1$
as $K\to\infty$. Thus, $(1-\varepsilon) \rme^{m/K} G(K)<1$ for $K$
sufficiently large, and we have proved
$\cPE^{\eta}[ \exp(\csigma_0/K) ]
\leq K \eta(V)$. To complete the proof, is suffices to note that
$\cPE^{\eta}[ \exp(\{\csigma_{1} - \csigma_0\}/K
)]
=\cPE^{\nu}[ \exp(\csigma_0/K) ]$ and
$\nu(V)<\infty$.
\end{pf}

With these preliminaries out of the way, we now prove Theorem
\ref{thm:exponential-inequality}.
\begin{pf*}{Proof Theorem \ref{thm:exponential-inequality}}
Define the sequence $(\xi_\ell)_{\ell\geq0}$ as in Lemma
\ref{lem:block-sums}. We begin by splitting the sum $S_n \eqdef
\sum_{i=1}^{n}\{f(X_i)-\pi(f)\}$ into three different terms:
%
%
\begin{eqnarray}\label{eq:dec1}
S_n &=&
\sum_{j=1}^{\csigma_0\wedge n-1} \{f(X_j) - \pi(f)\}
+ \sum_{k=0}^{i(n)-1}\xi_k\nonumber\\[-8pt]\\[-8pt]
&&{}+ \sum_{j=l(n)\wedge n}^{n} \{f(X_j) - \pi(f)\} ,\nonumber
\end{eqnarray}
where $i(n)\eqdef\sum_{k=1}^{\infty} \mathbh{1}_{\{ \csigma_k \leq
n \}}$
and $l(n)\eqdef\csigma_{i(n)}$. Using
(\ref{eq:bound-regeneration-time}), we have for $t>0$
%
%
\begin{eqnarray}
\label{eq:adamczak-1}
&&
\cPP^{\eta} \Biggl[ \Biggl|\sum_{j=1}^{\csigma_0\wedge n-1}
\{ f(X_j) - \pi(f) \} \Biggr| \geq t \Biggr]\nonumber\\
&&\qquad\le\cPP^{\eta} [ \csigma_0 \geq t / 2|f|_\infty]
\nonumber\\[-8pt]\\[-8pt]
&&\qquad\leq\cPE^{\eta} [ \exp(\csigma_0/K) ]
\exp\biggl( -\frac{t}{2K|f|_\infty}\biggr)\nonumber\\
&&\qquad\leq K \eta(V)
\exp\biggl( -\frac{t}{2K|f|_\infty}\biggr) .\nonumber
\end{eqnarray}
This bounds the first term of (\ref{eq:dec1}).
To bound the last term of (\ref{eq:dec1}), we proceed as in the
proof of \cite{adamczak2008}, Lemma 3. First note that, for any $t > 1$,
\begin{eqnarray*}
\cPP^{\eta}[ n - l(n)\wedge n +1 \geq t ] &=&
\cPP^{\eta} [ l(n) \leq n+1-t ] \\
&=&\sum_{\ell=0}^{n}
\cPP^{\eta}[ \csigma_\ell\leq n+1 - t, i(n) = \ell]\\
&=&
\sum_{\ell=0}^{n}
\cPP^{\eta}[ \csigma_\ell\leq n+1-t,
\csigma_{\ell+1} > n ] .
\end{eqnarray*}
Recall that the inter-regeneration time $\csigma_{\ell+1}- \csigma
_{\ell}$
is independent from $\csigma_0,\ldots,\csigma_\ell$, and
$(\csigma_{\ell+1}-\csigma_{\ell})_{\ell\ge0}$ are identically distributed
(see the proof of Lemma \ref{lem:block-sums}). Thus,
\begin{eqnarray*}
\cPP^{\eta}[ \csigma_\ell\leq n+1-t,
\csigma_{\ell+1} > n ] &=&
\sum_{k=0}^{\lfloor n+1-t \rfloor}
\cPP^{\eta} [ \csigma_\ell= k,
\csigma_{\ell+1}-\csigma_{\ell} > n -k] \\
&=& \sum_{k=0}^{\lfloor n+1-t \rfloor}
\cPP^{\eta} [ \csigma_\ell= k ]
\cPP^{\eta} [ \csigma_{1}-\csigma_{0} > n -k] .
\end{eqnarray*}
But as $\csigma_\ell< \csigma_{\ell+1}$ for all $\ell\geq0$,
we have $\sum_{\ell=0}^n \cPP^{\ceta}[
\csigma_\ell=k ] \leq1$ for all $k$, so that
\begin{eqnarray*}
\cPP^{\eta}[ n - l(n)\wedge n +1 \geq t ] &\leq&
\sum_{k=0}^{\lfloor n+1-t \rfloor}
\cPP^{\eta}[ \csigma_{1}-\csigma_{0} > n -k]
\\
&\leq& \sum_{k=\lceil t-1\rceil}^\infty
\cPP^{\eta}[ \csigma_{1}-\csigma_{0} \ge k]\\
&\leq&\cPE^{\eta}\bigl[ \rme^{\{\csigma_1-\csigma_0\}/K} \bigr]
\sum_{k=\lceil t-1\rceil}^\infty\rme^{-k/K}
\\
&\leq& \biggl(\frac{K \rme^{1/K}}{1-\rme^{-1/K}}\biggr)
\rme^{-t/K},
\end{eqnarray*}
where we have used (\ref{eq:bound-regeneration-time}).
We therefore find that for $t>2|f|_\infty$
%
%
\begin{eqnarray}\quad
\label{eq:adamczak-2}
\cPP_{\eta}\Biggl[ \Biggl| \sum_{j=l(n)\wedge n}^n
\{f(X_j) - \pi(f)\} \Biggr|
\geq t \Biggr]
&\leq&\cPP_{\eta} [ n - l(n)\wedge n + 1 \geq t /
2|f|_\infty] \nonumber\\[-8pt]\\[-8pt]
&\leq& K \exp\biggl( -\frac{t}{2K|f|_\infty}\biggr)\nonumber
\end{eqnarray}
(recall that the constant $K$ changes from line to line). But we may
clearly choose $K$ sufficiently large that $K\rme^{-1/K}\ge1$, so that
(\ref{eq:adamczak-2}) holds for any $t>0$.

It remains to bound the middle term in (\ref{eq:dec1}).
As $i(n)\le n$, we can estimate
\[
\Biggl|\sum_{k=0}^{i(n)-1} \xi_k\Biggr| \le
\max_{0\le j\le\lfloor n/2\rfloor}
\Biggl|\sum_{k=0}^{j} \xi_{2k}\Biggr| +
\max_{0\le j\le\lfloor n/2\rfloor}\Biggl|
\sum_{k=0}^{j} \xi_{2k+1}\Biggr|.
\]
Both terms on the right-hand side of this expression are identically
distributed. We can therefore estimate
using Etemadi's inequality (\cite{billingsley1987}, Theorem 22.5),
\[
\cPP^\eta\Biggl[ \Biggl|\sum_{k=0}^{i(n)-1} \xi_k\Biggr| \ge t
\Biggr] \le
8\max_{0\le j\le\lfloor n/2\rfloor}
\cPP^\eta\Biggl[ \Biggl|\sum_{k=0}^{j} \xi_{2k}\Biggr| \ge t/8
\Biggr].
\]
Note that $|\xi_k| \leq2|f|_\infty(\csigma_{k+1} - \csigma_k)$, so that
using (\ref{eq:bound-regeneration-time})
\[
(2K|f|_\infty)^2 \cPE^\eta\biggl(
\rme^{|\xi_k|/2K|f|_\infty} - 1- \frac{|\xi_k|}{2K|f|_\infty}
\biggr) \le
4K^3|f|_\infty^2.
\]
Using Bernstein's inequality (\cite{vandervaartwellner1996}, Lemma 2.2.11),
we obtain
\[
\cPP^\eta\Biggl[ \Biggl|\sum_{k=0}^{j} \xi_{2k}\Biggr| \ge t/8
\Biggr] \le
2 \exp\biggl(-\frac{1}{K}\frac{t^2}{(j+1)|f|_\infty^2
+t|f|_\infty}\biggr).
\]
We can therefore estimate for $t>0$
%
%
\begin{equation}
\label{eq:adamczak-3}
\cPP^\eta\Biggl[ \Biggl|\sum_{k=0}^{i(n)-1} \xi_k\Biggr| \ge t
\Biggr] \le
K \exp\biggl(-\frac{1}{K}\frac{t^2}{n|f|_\infty^2
+t|f|_\infty}\biggr).
\end{equation}
The proof is completed by combining (\ref{eq:adamczak-1}),
(\ref{eq:adamczak-2}) and (\ref{eq:adamczak-3}).
\end{pf*}

\subsection{\texorpdfstring{Proof of Theorem
\protect\ref{thm:azuma-hoeffding-v-unif}}{Proof of Theorem 14}}
\label{sec:proof-azuma-hoeffding-v-unif}

Assume without loss of generality that\break $\PEs_\theta[f({Y}_{0}^{s})]=0$.
To prove the result, it suffices to bound each term in the decomposition
\[
\sum_{i=1}^n f({Y}_{i}^{i+s}) =
\sum_{j=0}^s \Biggl( \sum_{i=1}^n \xi_{i,j} \Biggr)
+ \sum_{i=1}^n
\PE_\theta^\nu(f({Y}_{i}^{i+s})|X_0^{i-1},Y_0^{i-1}
),
\]
where we have defined for any $0\le j\le s$ and $i\geq1$
\[
\xi_{i,j} \eqdef
\PE_\theta^\nu(f({Y}_{i}^{i+s})
|X_0^{i+j},Y_0^{i+j}) -
\PE_\theta^\nu(f({Y}_{i}^{i+s})
|X_0^{i+j-1},Y_0^{i+j-1}).
\]
By construction, $(\xi_{i,j})_{1\le i\le n}$ are martingale increments
for each $j$, and $|\xi_{i,j}|_\infty\le2|f|_\infty$. Therefore,
by the Azuma--Hoeffding inequality (\cite{williams1991}, page 237), we have
\[
\PP_\theta^\nu\Biggl(\Biggl|\sum_{i=1}^n \xi_{i,j}\Biggr|
\ge t\Biggr) \le
2 \exp\biggl(
-\frac{t^2}{8n|f|_\infty^2}
\biggr)
\]
for each $0\le j\le s$. On the other hand, note that
$\PE_\theta^\nu(f({Y}_{i}^{i+s})|X_0^{i-1},Y_0^{i-1}
)=F(X_{i-1})$ for all $i$, where $F$ satisfies
$\pi_\theta(F)=0$ (as we assumed $\PEs_\theta[f({Y}_{0}^{s})]=0$)
and $|F|_\infty\le|f|_\infty$. The result therefore follows by applying
Theorem \ref{thm:exponential-inequality}.


%
\printaddresses

\end{document}